\documentclass[12pt,smallextended]{svjour3}
\smartqed
\usepackage{lineno}

\journalname{Journal of Scientific Computing}
\usepackage{amssymb,amsmath,amsfonts}
\usepackage{hyperref}
\hypersetup{
    colorlinks=true,
    linkcolor=blue,
    filecolor=blue,
    urlcolor=blue,
    citecolor=cyan,
    pdfstartview=FitW,
    bookmarks=true,
}
\usepackage{array}
\usepackage{epsfig}
\usepackage{mathptmx}
\usepackage{diagbox}
\usepackage{graphicx}
\usepackage{color}
\usepackage{caption}
\usepackage{wrapfig}
\usepackage{graphics}
\usepackage{enumerate}
\usepackage{amsxtra}
\usepackage{latexsym}
\usepackage{multirow}
\usepackage{slashbox}
\usepackage{subfigure}
\usepackage{epstopdf}
\usepackage{pdfpages}
\usepackage{algorithm}
\usepackage{algorithmic}
\usepackage{url}
\newcommand{\thmlist}{
\begin{list}{Step 1}
{\setlength{\leftmargin}{0.6 in}\setlength{\labelwidth} {0.5 in}}}

\newcommand{\alglist}{
\begin{list}{Step 1}
{\setlength{\leftmargin}{1.1 in} \setlength{\labelwidth}{1.0 in}}}

 \renewcommand{\proof} {\noindent {\bf Proof.} \quad}

\renewcommand{\subtitle}[1]{\color{blue}}

\def\red#1{\color{red}{#1}\color{black}}

\begin{document}



\title{Generalized continuation Newton methods and the trust-region updating strategy
for the underdetermined system}
\titlerunning{Generalized continuation Newton methods and the trust-region
updating strategy}
\author{Xin-long Luo \textsuperscript{$\ast$} \and Hang Xiao}
\authorrunning{Luo, Xiao}

\institute{
     Xin-long Luo, Corresponding author
     \at
     School of Artificial Intelligence,
     Beijing University of Posts and Telecommunications, P. O. Box 101,
     Xitucheng Road  No. 10, Haidian District, 100876, Beijing China \\
     \email{luoxinlong@bupt.edu.cn}
     \and
     Hang Xiao
     \at
     School of Artificial Intelligence,
     Beijing University of Posts and Telecommunications, P. O. Box 101,
     Xitucheng Road  No. 10, Haidian District, 100876, Beijing China \\
     \email{xiaohang0210@bupt.edu.cn}
}

\date{Received: date / Accepted: date}

\maketitle

\begin{abstract}
This paper considers the generalized continuation Newton method and the
trust-region updating strategy for the underdetermined system of nonlinear 
equations. Moreover, in order to improve its computational efficiency, the new 
method will not update the Jacobian matrix when the current Jacobian matrix 
performs well. The numerical results show that the new method is more robust 
and faster than the traditional optimization method such as the 
Levenberg-Marquardt method (a variant of trust-region methods, the built-in 
subroutine fsolve.m of the MATLAB R2020a environment). The computational time of 
the new method is about 1/8 to 1/50 of that of fsolve. Furthermore, it also 
proves the global convergence and the local superlinear convergence of the new 
method under some standard assumptions.
\end{abstract}




\section{Introduction}
In engineering fields, we often need to solve the underdetermined system of equations
as follows:
\begin{align}
      F(x) = 0, \label{UNLEFX}
\end{align}
where $F: \; \Re^{n} \to \Re^{m}$ and $m < n$. For example, this problem arises from
finding the initial feasible point of the following differential-algebraic equations
\cite{AP1998,BCP1996,HW1996,LF2000,LL2001}:
\begin{align}
  & \frac{dx}{dt} = h(x, \, y),   \label{ODE} \\
  & g(x) = 0.    \label{ALBEQ}
\end{align}
Another case comes from the feasible direction method for solving the following
nonlinearly constrained optimization
problem \cite{NW1999,SY2006}
\begin{align}
    \min_{x \Re^{n}}  \;  r(x) \; \text{subject to} \; c(x) = 0, \label{NCEOPT}
\end{align}
where $r: \; \Re^{n} \to \Re$ and $c: \Re^{n} \to \Re^{m}, \; m < n$.

\vskip 2mm

The main difficulty of the undertermined system is the singularity 
$J(x)^{T}J(x)$ \cite{FY2005,Griewank1985,QACT2017,YF2001}, where $J = F'$ is the
Jacobian function of $F$. When $m = n$ and the Jacobian matrix 
$J(x)$ is nonsingular, there are many popular traditional optimization
methods \cite{CGT2000,DS2009,Higham1999,Kelley2018,NW1999,Yuan1998}
and the classical homotopy continuation methods \cite{AG2003,Doedel2007,OR2000,WSMMW1997}
to solve it.

\vskip 2mm

For the traditional optimization methods such as the trust-region methods
and the line search methods, the solution $x^{\ast}$ of the nonlinear system
\eqref{UNLEFX} is found via solving the following equivalent nonlinear least-squares
problem
\begin{align}
    \min_{x \in \Re^{n}} f(x) = \frac{1}{2}\|F(x)\|^2,         \label{UOPTF}
\end{align}
where $\|\cdot\|$ denotes the Euclidean vector norm or its induced
matrix norm throughout this paper. Generally speaking, the traditional optimization
methods based on the merit function \eqref{UOPTF} are efficient for the large-scale
problems when $J(x_{k})^{T}J(x_k) \, (k = 0, \, 1, \, \ldots)$ are nonsingular, 
since they have the local superlinear convergence near the solution $x^{\ast}$
\cite{CGT2000,NW1999}.

\vskip 2mm

However, the line search method based on the classical Gauss-Newton method 
will confront some problems when $J(x_k)^{T}J(x_k)$ is singular, since it obtains 
the search direction $d_{k}$ by solving the following linear equations: 
\begin{align}
        J(x_k)^{T}J(x_k)d_{k} = -J(x_k)^{T}F(x_k).  \nonumber 
\end{align}
Furthermore, the termination condition
\begin{align}
   \|\nabla f(x_k)\| = \|J(x_k)^{T} F(x_k)\| < \epsilon, \label{TUNOPT}
\end{align}
may lead the methods based on the merit function \eqref{UOPTF} to early stop 
far away from the solution $x^{\ast}$. This can be illustrated as follows. 
We consider
\begin{align}
  F(x) = Ax = 0,  \hskip 2mm  A =
     \begin{bmatrix}
     1 & 0 \\
     0 & 10^{-6}
    \end{bmatrix}. \label{TINYEXA}
\end{align}
It is not difficult to know that the linear system \eqref{TINYEXA} has  a unique
solution $x^{\ast} = (0, \, 0)$. If we set $\epsilon = 10^{-6}$,
the traditional optimization methods will early stop far away from
$x^{\ast}$ provided that $x_{k} = (0, \; c), \; c < 10^{6}$.

\vskip 2mm

For the classical homotopy methods, the solution $x^{\ast}$ of the nonlinear
system \eqref{UNLEFX} is found via constructing the following homotopy function
\begin{align}
  H(x, \, \lambda) = (1- \lambda)G(x) + \lambda F(x), \label{HOMOFUN}
\end{align}
and attempting to trace an implicitly defined curve $\lambda(t) \in H^{-1}(0)$
from the starting point $(x_{0}, \, 0)$ to a solution $(x^{\ast}, \, 1)$ by the
predictor-corrector methods \cite{AG2003,Doedel2007}, where the zero point of
the artificial smooth function $G(x)$ is known. Generally speaking, the homotopy
continuation methods are more reliable than the merit-function methods and they
are very popular in engineering fields \cite{Liao2012}. The disadvantage of the
classical homotopy methods is that they require significantly more function and
derivative evaluations, and linear algebra operations than the merit-function
methods since they need to solve many auxiliary nonlinear systems during the
intermediate continuation process.

\vskip 2mm

In order to overcome this shortcoming of the traditional homotopy methods, we
consider the special continuation method based on the following generalized 
Newton flow \cite{AS2015,Branin1972,Davidenko1953,LXL2020,Tanabe1979}
\begin{align}
   \frac{dx(t)}{dt} = - J(x)^{+}F(x), \;  x(t_0) = x_0,  \label{NEWTONFLOW}
\end{align}
where $J(x)^{+}$ is the Moore-Penrose generalized inverse of the Jacobian
matrix $J(x)$ (p. 11, \cite{SY2006} or p. 290, \cite{GV2013}). Then, we construct
a special ODE method with the adaptively time-stepping scheme based on the
trust-region updating strategy to trace the trajectory of the generalized 
Newton flow \eqref{NEWTONFLOW}. Consequently, we obtain a solution $x^{\ast}$ of the
underdetermined nonlinear system \eqref{UNLEFX}.

\vskip 2mm

The rest of this article is organized as follows. In the next section, we
consider the generalized continuation Newton method with the adaptively time-stepping
scheme and the updating technique of the Jacobian matrix based on the trust-region updating
strategy for the underdetermined system of nonlinear equations. In section 3,
under the standard assumptions, we prove the global convergence and the local
superlinear convergence of the new method. In section 4, some promising numerical
results of the new method are also reported, in comparison to the Levenberg-Marquardt
method (a variant of the trust-region methods,
the built-in subroutine fsolve.m of the MATLAB R2020a environment)
 \cite{FY2005,Levenberg1944,MATLAB,Marquardt1963,More1978,YF2001}).
Finally, some conclusions and the discussions are given in section 5. Throughout
this article, we assume that $F(\cdot)$ exists the zero point $x^{\ast}$.

\vskip 2mm

\section{Continuation Newton methods}

In this section, based on the trust-region updating strategy, we construct an
adaptively time-stepping scheme for the continuation Newton method to trace the
trajectory of the generalized Newton flow and obtain its equilibrium point $x^{\ast}$.

\subsection{The generalized continuous Newton flow} \label{SUBSECNF}

If we consider the damped Newton method with the line search strategy for
the nonlinear system \eqref{UNLEFX} \cite{Kelley2003,NW1999}, we have
\begin{align}
   x_{k+1} = x_{k} - \alpha_{k} J(x_{k})^{+} F(x_{k}). \label{NEWTON}
\end{align}
We denote $o(\alpha)$ as the higher-order infinitesimal of $\alpha$, that is to
say,
\begin{align}
   \lim_{\alpha \to 0} \frac{o(\alpha)}{\alpha} = 0. \nonumber
\end{align}
In equation \eqref{NEWTON}, if we regard  $x_{k} = x(t_{k})$ and
$x_{k+1} = x(t_{k} + \alpha_{k}) + o(\alpha_{k})$, we obtain the continuous
Newton flow \eqref{NEWTONFLOW} when $\alpha_{k} \to 0$. Actually, if we apply
an iteration with the explicit Euler method \cite{SGT2003} to the
generalized Newton flow \eqref{NEWTONFLOW}, we also obtain the damped Newton
method \eqref{NEWTON}. Since the rank of the Jacobian matrix $J(x)$ may be not
full, we reformulate the generalized Newton flow \eqref{NEWTONFLOW} as the more
general formula:
\begin{align}
  -J(x)\frac{dx(t)}{dt} = F(x), \hskip 2mm  x(t_0) = x_0. \label{DAEFLOW}
\end{align}

\vskip 2mm

The continuous Newton flow \eqref{DAEFLOW} is an old method and can be backtracked
to Davidenko's work \cite{Davidenko1953} in 1953. After that, it was investigated
by Branin \cite{Branin1972}, Deuflhard et al \cite{PHP1975}, Tanabe \cite{Tanabe1979}
and Kalaba et al \cite {KZHH1977} in 1970s, and applied to nonlinear boundary
problems by Axelsson and Sysala \cite{AS2015} recently. The continuous and even
growing interest in this method originates from its some nice properties. One of
them is that the solution $x(t)$ of the continuous Newton flow converges to
the steady-state solution $x^{\ast}$ from any initial point $x_{0}$, as described
by the following property \ref{PRODAEFLOW}.

\vskip 2mm

\begin{property} (Branin \cite{Branin1972} and Tanabe \cite{Tanabe1979})
\label{PRODAEFLOW} Assume that $x(t)$ is the solution of the continuous Newton
flow \eqref{DAEFLOW}, then $f(x(t)) = \|F(x)\|^{2}$ converges to zero when
$t \to \infty$. That is to say, for every limit point $x^{\ast}$ of $x(t)$, it is
also a solution of the underdetermined system \eqref{UNLEFX}. Furthermore, every
element $F^{i}(x)$ of $F(x)$ has the same convergence rate $e^{-t}$
and $x(t)$ can not converge to the solution $x^{\ast}$ of the underdetermined system
\eqref{UNLEFX} on the finite interval when the initial point $x_{0}$ is not a
solution of the underdetermined system \eqref{UNLEFX}.
\end{property}
\proof Assume that $x(t)$ is the solution of the continuous Newton flow
\eqref{DAEFLOW}, then we have
\begin{align}
    \frac{d}{dt} \left(e^{t}F(x)\right) = e^{t} J(x) \frac{dx(t)}{dt}
    + e^{t} F(x) = 0. \nonumber
\end{align}
Consequently, we obtain
\begin{align}
     F(x(t)) = F(x_0)e^{-t}. \label{FUNPAR}
\end{align}
From equation \eqref{FUNPAR}, it is not difficult to know that every element
$F^{i}(x)$ of $F(x)$ converges to zero with the linear convergence
rate $e^{-t}$ when $t \to \infty$. Thus, if the solution $x(t)$ of the continuous
Newton flow \eqref{DAEFLOW} belongs to a compact set, it has a limit point
$x^{\ast}$ when $t \to \infty$, and this limit point $x^{\ast}$ is
also a solution of the underdetermined system \eqref{UNLEFX}.

\vskip 2mm

If we assume that the solution $x(t)$ of the continuous Newton flow
\eqref{DAEFLOW} converges to the solution $x^{\ast}$ of the underdetermined system
\eqref{UNLEFX} on the finite interval $(0, \, T]$, from equation \eqref{FUNPAR},
we have
\begin{align}
     F(x^{\ast}) = F(x_{0}) e^{-T}. \label{FLIMT}
\end{align}
Since $x^{\ast}$ is a solution of the underdetermined system \eqref{UNLEFX}, we have
$F(x^{\ast}) = 0$. By substituting it into equation \eqref{FLIMT}, we
obtain
\begin{align}
     F(x_{0}) = 0. \nonumber
\end{align}
Thus, it contradicts the assumption that $x_{0}$ is not a solution of the underdetermined
system \eqref{UNLEFX}. Consequently, the solution $x(t)$ of the continuous Newton flow
\eqref{DAEFLOW} can not converge to the solution $x^{\ast}$ of the underdetermined system
\eqref{UNLEFX} on the finite interval. \qed

\vskip 2mm

\begin{remark}
The inverse $J(x)^{-1}$ of the Jacobian matrix $J(x)$ can be regarded as the
preconditioner of $F(x)$ such that the solution elements $x_{i}(t)\, (i = 1, \, 2, \ldots, n)$
of the continuous Newton flow \eqref{NEWTONFLOW} have the roughly same convergence
rates and it mitigates the stiff property of the ODE \eqref{NEWTONFLOW}
(the definition of the stiff problem can be found in \cite{HW1996} and references
therein). This property is very useful since it makes us adopt the explicit ODE
method to trace the trajectory of the Newton flow.
\end{remark}

\vskip 2mm

Actually, if we consider $F(x) = Ax$, from the ODE \eqref{DAEFLOW},
we have
\begin{align}
  A \frac{dx}{dt} = - A x, \; x(0) = x_{0}.  \label{LINAPP}
\end{align}
By integrating the linear ODE \eqref{LINAPP}, we obtain
\begin{align}
   x(t) = e^{-t} x_{0}.       \label{LINSOL}
\end{align}
From equation \eqref{LINSOL}, we know that the solution $x(t)$ of the ODE
\eqref{LINAPP} converges to zero exponentially with the same rate $e^{-t}$
when $t$ tends to infinity.

\subsection{The generalized continuation Newton method}

From subsection \ref{SUBSECNF}, we know that the solution $x(t)$ of the
generalized continuous Newton flow \eqref{DAEFLOW} has the nice global convergence
property. On the other hand, when the Jacobian matrix $J(x)$ is singular or nearly
singular, the ODE \eqref{DAEFLOW} is the system of differential-algebraic equations
(DAEs) and its trajectory can not be efficiently followed by the general ODE method
such as the backward differentiation formulas (the built-in subroutine ode15s.m of the
MATLAB environment \cite{AP1998,BCP1996,HW1996,MATLAB,SGT2003}). Thus, we need to
construct the special method to handle this problem.

\vskip 2mm

Since the continuous Newton flow \eqref{DAEFLOW} is intrinsically
a nonlinear diminishing system for the energy function $f(x(t)) = \|F(x(t))\|^{2}$,
it can be integrated by the strong stability preserving methods \cite{GS1998,GST2001}
and the steady-state solution $x^{\ast}$ can be obtained after the long time
integration. Here, we consider another approach based on the traditional optimization
methods for problem \eqref{DAEFLOW}. We expect that the new method has the global
convergence as the homotopy continuation method and the fast convergence rate near
the steady-state solution $x^{\ast}$ as the merit-function method. In order to
achieve these two aims, we construct the special continuation Newton method with the
new step size $\alpha_{k} = \Delta t_{k}/(1+\Delta t_{k})$ and the time step size
$\Delta t_{k}$ is adaptively adjusted by the trust-region updating strategy for
problem \eqref{DAEFLOW}.

\vskip 2mm

Firstly, we apply the implicit Euler method to the continuous Newton flow
\eqref{DAEFLOW} \cite{AP1998,BCP1996}, then we obtain
\begin{align}
  J(x_{k+1}) \frac{x_{k+1}-x_{k}}{\Delta t_k} & =  -F(x_{k+1}).  \label{IMEDAE}
\end{align}
The scheme \eqref{IMEDAE} is an implicit method and it needs to solve a system
of nonlinear equations at every iteration. To avoid solving the system of nonlinear
equations, we replace $J(x_{k+1})$ with $J(x_{k})$ and substitute $F(x_{k+1})$ with
its linear approximation $F(x_k)+J(x_k)(x_{k+1}-x_{k})$ in equation \eqref{IMEDAE}.
Thus, we obtain the generalized continuation Newton method as follows:
\begin{align}
       J(x_{k})s_{k}^{N} = -  F(x_k),
      \; x_{k+1} = x_{k} + \frac{\Delta t_k}{1+\Delta t_k}s_{k}^{N}. \label{GCNM}
\end{align}

\vskip 2mm

The linear system \eqref{GCNM} is underdetermined. That is to say, the row rank
of $J(x_{k})$ is less than the number of the variable $s_{k}^{N}$. Thus, the linear
system \eqref{GCNM} may have many solutions or no solution. For simplicity,
we assume that the Jacobian matrix $J(x_{k})$ is full row rank. That is to say,
the row rank of $J(x_{k})$ equals $m$. In order to obtain the nearest point
$x_{k+1}$ of $x_{k}$ under the constraint \eqref{GCNM}, we solve the following
shortest distance problem:
\begin{align}
      \min_{s^{N} \in \Re^{n}} \; \left\|s^{N}\right\|^{2}, \; \text{subject to} \;
      J_{k}s^{N} = - F_{k}, \label{CSDP}
\end{align}
where $J_{k}$ equals $J(x_k)$ or its approximation and $F_{k} = F(x_k)$. By
using the Lagrangian multiplier method \cite{SY2006}, it is not difficult to
obtain the solution $s_{k}^{N}$ of problem \eqref{CSDP} as follows:
\begin{align}
     s_{k}^{N} = - J_{k}^{+} F_{k}, \;
     J_{k}^{+} = J_{k}^{T} \left(J_{k}J_{k}^{T}\right)^{-1}, \label{SMEDAE}
\end{align}
where $J_{k}^{+}$ is the pseudo-inverse of $J_{k}$. Thus, from equations
\eqref{GCNM} and \eqref{SMEDAE}, we obtain the generalized continuation Newton
method for the underdetermined system \eqref{UNLEFX} as follows:
\begin{align}
     x_{k+1} = x_{k} - \frac{\Delta t_k}{1+\Delta t_k} J_{k}^{+} F_{k}.
     \label{GICNM}
\end{align}

\vskip 2mm

The matrix $J_{k}J_{k}^{T}$ may be ill-conditioned. Thus, the Cholesky
factorization method may fail to solve the linear system \eqref{SMEDAE} for the
large-scale problem. Therefore, we use the QR decomposition
(pp. 247-248, \cite{GV2013}) to solve it as follows:
\begin{align}
     J_{k}^{T} = Q_{k}R_{k},  \;  R_{k}^{T} d_{k} = - F_{k},  \;
     s_{k}^{N} = Q_{k} d_{k}, \;
     s_{k} = \frac{\Delta t_{k}}{1 + \Delta t_{k}} s_{k}^{N},  \label{SKJTQR}
\end{align}
where $Q_{k} \in \Re^{n \times m}$ satisfies $Q_{k}^{T}Q_{k} = I$ and
$R_{k} \in \Re^{m \times m}$ is an upper triangle matrix.

\vskip 2mm

\begin{remark}
The generalized continuation Newton method \eqref{GICNM} is similar to the
damped Newton method \eqref{NEWTON} if we let $\alpha_{k} =
\Delta t_k/(1+\Delta t_k)$ in equation \eqref{GICNM}. However, from the view of
the ODE method, they are different. The damped Newton method \eqref{NEWTON} is
obtained by the explicit Euler scheme applied to the generalized continuous Newton
flow \eqref{DAEFLOW}, and its time step size $\alpha_k$ is restricted by the
numerical stability \cite{HW1996,SGT2003}. That is to say, for the linear
test equation $dx/dt = - \lambda x$, its time step size $\alpha_{k}$ is
restricted by the stable region $|1-\lambda \alpha_{k}| \le 1$. Therefore, the
large time step can not be adopted in the steady-state
phase. The generalized continuation Newton method \eqref{GICNM} is obtained by
the implicit Euler method and its linear approximation applied to the continuous Newton
flow \eqref{DAEFLOW}, and its time step size $\Delta t_k$ is not restricted
by the numerical stability. Therefore, the large time step can be
adopted in the steady-state phase and it mimics the generalized Newton method near the solution
$x^{\ast}$ such that it has the fast local convergence rate. The most
of all, $\alpha_{k} = \Delta t_{k}/(\Delta t_{k} + 1)$ in equation \eqref{GICNM} is
favourable to adopt the trust-region updating strategy for adaptively adjusting the
time step size $\Delta t_{k}$ such that the generalized continuation
Newton method \eqref{GICNM} accurately traces the trajectory of the generalized
continuous Newton flow in the transient-state phase and achieves the fast convergence
rate near the equilibrium point $x^{\ast}$.
\end{remark}

\vskip 2mm

\begin{remark}
We denote $r(x) = \|F(x)\|$ and the generalized Newton direction
$s_{k}^{N}$ as
\begin{align}
      s_{k}^{N} = - J(x_{k})^{+}F(x_{k}). \label{GND}
\end{align}
Then, when $F(x_{k}) \neq 0$, we have
\begin{align}
      \nabla r(x_{k})^{T} s_{k}^{N}
      = - \frac{F(x_{k})^{T}J(x_{k})}{\|F(x_{k})\|}
      \left(J(x_{k})^{+}F(x_{k})\right)
      = - \|F(x_{k})\| < 0.\label{DESD}
\end{align}
That is to say, the generalized Newton direction $s_{k}^{N}$ is a descent
direction of $r(x_{k})$.
\end{remark}

\vskip 2mm

\subsection{The trust-region updating strategy}

Another issue is how to adaptively adjust the time step size $\Delta t_k$
at every iteration. There is a popular way to control the time step size 
based on the trust-region updating strategy
\cite{CGT2000,Deuflhard2004,Higham1999,Luo2009,Luo2010,LLT2007,LLS2021,LY2021,Yuan2015}.
Its main idea is that the time step size $\Delta t_{k+1}$ will be enlarged
when the linear model $F(x_{k}) + J_{k}s_{k}$ approximates $F(x_{k}+s_{k})$ well, and
$\Delta t_{k+1}$ will be reduced when $F(x_k) + J_{k}s_{k}$ approximates
$F(x_{k}+s_{k})$ badly.

\vskip 2mm 

In practice, we enlarge or reduce the time step size $\Delta t_k$ at every 
iteration according to the following ratio:
\begin{align}
  \rho_k = \frac{\|F(x_{k})\|-\|F(x_{k}+s_{k})\|}
   {\|F(x_{k})\| - \|F(x_{k})+ J_{k}s_{k}\|}.  \label{RHOKO}
\end{align}
From the computational formula \eqref{GICNM} of the search step $s_{k}$, we can 
save the computational time of the predicted model $F(x_k) + J_{k}s_{k}$ by the 
following simplified formula:   
\begin{align}
   F(x_k) + J_{k}s_{k} = F(x_k) - \frac{\Delta t_{k}}{1+\Delta t_{k}} F(x_k)
   = \frac{1}{1+\Delta t_{k}} F(x_k). \label{PREMFK1}
\end{align}
Thus, from equations \eqref{RHOKO}-\eqref{PREMFK1}, we rewrite the computational 
formula \eqref{RHOKO} of $\rho_{k}$ as 
\begin{align}
  \rho_k = \frac{\|F(x_{k})\|-\|F(x_{k}+s_{k})\|}
   {(\Delta t_{k}/(1+\Delta t_{k}))\|F(x_k)\|}.  \label{RHOK}
\end{align}
Therefore, according to the computation formula \eqref{RHOK} of $\rho_{k}$ 
between the actual reduction and the predicted reduction, we give 
a particular adjustment strategy of $\Delta t_{k}$ as follows:
\begin{align}
   \Delta t_{k+1} =
     \begin{cases}
    \gamma_1 \Delta t_k, &{\text{if} \;  \left|1- \rho_k \right| \le \eta_1,}\\
    \Delta t_k, & {\text{else if} \; \eta_1 < \left|1 - \rho_k \right| < \eta_2,}\\
    \gamma_2 \Delta t_k, &{\text{others},}
    \end{cases} \label{TSK1}
\end{align}
where the constants are selected as $\gamma_{1} = 2, \; \gamma_{2} = 0.5, \;
\eta_{1} = 0.25, \; \eta_{2} = 0.75$, according to our numerical experiments.

\vskip 2mm

\begin{remark}
This new time-stepping scheme based on the trust-region updating
strategy has some advantages compared to the traditional line search strategy
\cite{Luo2005}. If we use the line search strategy and the damped Newton method
\eqref{NEWTON} to track the trajectory $z(t)$ of the generalized continuous
Newton flow \eqref{DAEFLOW}, in order to achieve the fast convergence rate in
the steady-state phase, the time step size $\alpha_{k}$ of the damped Newton
method is tried from 1 and reduced by half with many times at every iteration.
Since the linear model $F(x_{k}) + J(x_{k})s_{k}$ may not approximate
$F(x_{k}+s_{k})$ well in the transient-state phase, the time step size $\alpha_{k}$
will be small. Consequently, the line search strategy consumes the unnecessary
trial steps in the transient-state phase. However, the selection scheme of the
time step size based on the trust-region updating strategy \eqref{RHOK}-\eqref{TSK1}
can overcome this shortcoming.
\end{remark}

\vskip 2mm

\subsection{The updating technique of the Jacobian matrix}

\vskip 2mm

For a system of nonlinear equations, the computational time of
the Jacobian matrix is heavy if we update the Jacobian matrix $J(x_{k})$
at every iteration. In order to save the computational time of
the Jacobian evaluation, similarly to the switching preconditioned technique
\cite{LXLZ2021}, we set $J_{k+1} = J_{k}$ when $F_{k} + J_{k}s_{k}$
approximates $F(x_{k}+s_{k})$ well. Otherwise, we update $J_{k+1} = J(x_{k+1})$.
An effective updating strategy is give by
\begin{align}
     J_{k+1} = \begin{cases}
                 J_{k},  \; \text{if} \; |1- \rho_{k}| \le \eta_{1}, \\
                 J(x_{k+1}), \; \text{otherwise},
               \end{cases} \label{UPDJK1}
\end{align}
where $\rho_{k}$ is defined by equation \eqref{RHOK} and $\eta_{1} = 0.25$. 
In practice, in order to save the computational time of decomposing the matrix 
$J_{k}$ when $J_{k-1}$ performs well, i.e. $|1-\rho_{k}| \le \eta_{1}$, 
according to the updating formula \eqref{UPDJK1}, we set $R_{k} = R_{k-1}$ and 
$Q_{k} = Q_{k-1}$ in equation \eqref{SKJTQR}. 

\vskip 2mm

For a real-world problem, the analytical Jacobian $J(x_{k})$ may not be
offered. Thus, in practice, we replace the Jacobian matrix $J(x_{k})$
with its difference approximation as follows:
\begin{align}
     J(x_{k}) \approx
    \left[\frac{F(x_{k} + \epsilon e_{1}) - F(x_{k})}{\epsilon}, \,
    \ldots, \, \frac{F(x_{k} + \epsilon e_{n}) - F(x_{k})}{\epsilon}\right],
    \label{NUMHESS}
\end{align}
where $e_{i}$ represents the unit vector whose elements equal zeros except
for the $i$-th element which equals 1, and the parameter $\epsilon$ can be selected as
$10^{-6}$ according to our numerical experiments.

\vskip 2mm

According to the above discussions, we give the detailed implementation of
the generalized continuation Newton method with the trust-region updating
strategy for the underdetermined system of nonlinear equations in Algorithm
\ref{ALGGCNMTR}.

\vskip 2mm

\begin{algorithm}
   \renewcommand{\algorithmicrequire}{\textbf{Input:}}
   \renewcommand{\algorithmicensure}{\textbf{Output:}}
   \caption{Generalized continuation Newton methods and the trust-region
   updating strategy for the underdetermined system  (The GCNMTr method)}
   \label{ALGGCNMTR}
   \begin{algorithmic}[1]
      \REQUIRE ~~ \\
      Function $F: \; \Re^{n} \to \Re^{m}, \; m \le n$,
      the initial point  $x_0$ (optional), and the tolerance $\epsilon$ (optional).
	  \ENSURE ~~ \\
      An approximation solution $x^{\ast}$ of nonlinear equations.
      \STATE Set the default $x_0 = \text{ones} (n, \, 1)$  and
      $\epsilon = 10^{-6}$ when $x_0$ or $\epsilon$ is not provided by the calling subroutine.
      \STATE Initialize the parameters: $\eta_{a} = 10^{-6}, \; \eta_1 = 0.25, \;
      \gamma_1 =2, \; \eta_2 = 0.75, \; \gamma_2 = 0.5, \; \text{maxit} = 400$.
      \STATE Set $\Delta t_0 = 10^{-2}$, flag\_success\_trialstep = 1, $\text{itc} = 0, \; k = 0$.
      \STATE Evaluate $F_{k} = F(x_{k})$ and $J_k = J(x_{k})$.
      \STATE By using the qr decomposition $[Q_{k},R_{k}] = \text{qr}(J_{k}^{T})$ of $J_{k}^{T}$,
      we obtain the orthogonal matrix $Q_{k}$ and the upper triangle matrix $R_{k}$.
      \STATE Set $\rho_{k}  = 1$.
      \WHILE{(itc $<$ maxit)}
         \IF{(flag\_success\_trialstep == 1)}
              \STATE Set itc = itc + 1.
              \STATE Compute $\text{Res}_{k} = \|F_{k}\|_{\infty}$.
              \IF{($\text{Res}_{k} < \epsilon$)}
                 \STATE break;
              \ENDIF
              \IF{$(|1-\rho_{k}| > 0.25)$}
                 \STATE Evaluate $J_k = J(x_{k})$.
                 \STATE By using the qr decomposition $[Q_{k},R_{k}] = \text{qr}(J_{k}^{T})$ of $J_{k}^{T}$,
                  we obtain the orthogonal matrix $Q_{k}$ and the upper triangle matrix $R_{k}$.
              \ENDIF
              \STATE By solving $R_{k}^{T}b_{k} = - F_{k}$ and $s_{k}^{N} = Q_{k}b_{k}$,
              we obtain the Newton step $s_{k}^{N}$.
         \ENDIF
         \STATE Set $s_{k} = {\Delta t_{k}}/{(1+\Delta t_{k})} \, s_{k}^{N}, \;
         x_{k+1} = x_k + s_k$.
         \STATE Evaluate $F(x_{k+1})$.
         \IF {$\left(\|F_{k}\| < \|F_{k} + J_{k}s_k\|\right)$}
           \STATE $\rho_{k} = -1$;
         \ELSE
           \STATE Compute the ratio $\rho_{k}$ from equation \eqref{RHOK}.
         \ENDIF
         \STATE Adjust the time step size $\Delta t_{k+1}$ according to the
         trust-region updating strategy \eqref{TSK1}.
         \IF{$(\rho_{k} \ge \eta_{a})$}
               \STATE Accept the trial point $x_{k+1}$. Set flag\_success\_trialstep = 1.
           \ELSE
               \STATE Set $x_{k+1}  = x_{k}$, $F_{k+1} = F_{k}$, $s_{k+1}^{N} = s_{k}^{N}$, 
               flag\_success\_trialstep = 0.
           \ENDIF
           \STATE Set $\rho_{k+1} = \rho_{k}$, $R_{k+1} = R_{k}$, $Q_{k+1} = Q_{k}$.
         \STATE Set $k \longleftarrow k+1$.
      \ENDWHILE
   \end{algorithmic}
\end{algorithm}

\section{Algorithm analysis}

In this section, we discuss some theoretical properties of Algorithm \ref{ALGGCNMTR}.
Firstly, we estimate the lower bound of the predicted reduction
$\|F(x_k)\|-\|F(x_{k})+ J_{k}s_{k}\|$, which is similar to that
of the trust-region method for the unconstrained optimization problem
\cite{Powell1975}.

\vskip 2mm

According to the theorem of the singular value decomposition
(pp. 76, \cite{GV2013}), for the matrix $J_{k} \in \Re^{m \times n}$,
there exist orthogonal matrices $U_{k} \in \Re^{m \times m}$ and
$V_{k} \in \Re^{n \times n}$ such that
\begin{align}
   U_{k}^{T}J_{k}V_{k} = \Sigma_{k} =
   \text{diag}\left(\sigma_{k}^{1}, \, \ldots, \, \sigma_{k}^{m}\right)
   \in \Re^{m \times n},    \label{SVD}
\end{align}
where $\sigma_{k}^{1} \ge \sigma_{k}^{2} \ge \cdots \ge \sigma_{k}^{m} \ge 0$.

\vskip 2mm

\begin{lemma} \label{LEMESTPRED}
Assume that it exists a positive constant $c_{\sigma}$ such that
\begin{align}
  \sigma_{k}^{min} \ge c_{\sigma}  \label{UBINVJ}
\end{align}
holds for all $k = 0, \, 1, \ldots$, where $\sigma_{k}^{min} = \sigma_{k}^{m}$
is the smallest singular value of $J_{k} \in \Re^{m \times n}$. Furthermore,
we suppose that $s_{k}$ is the solution of the generalized continuation Newton
method \eqref{SMEDAE}-\eqref{GICNM}. Then, we have the following estimation
\begin{align}
   \|F(x_k)\| - \|F(x_k) + J_{k}s_{k}\| =
    \frac{\Delta t_k}{1+\Delta t_k}\|F(x_k)\|. \label{PRERED}
\end{align}
\end{lemma}
\proof From equations \eqref{SVD}-\eqref{UBINVJ}, we have
\begin{align}
  J_{k}^{+}  =  J_{k}^{T}\left(J_{k}J_{k}^{T}\right)^{-1}
  = V_{k}^{T} \Sigma_{k}^{-1}U_{k}, \;
  \Sigma_{k}^{-1} = \text{diag}\left({1}/{\sigma_{k}^{1}}, \, \ldots, \,
  {1}/{\sigma_{k}^{m}}\right) \in \Re^{n \times m}. \label{GINVJK}
\end{align}
Thus, from equations \eqref{GICNM}, \eqref{SVD} and \eqref{GINVJK}, we have
\begin{align}
   \|F(x_k)+J_{k}s_{k}\|   = \|J_{k}s_{k} + F_{k}\|
   = \left\|- \frac{\Delta t_{k}}{1+\Delta t_{k}}J_{k}J_{k}^{+}F_{k}
   + F_{k} \right\| =  \frac{1}{1 + \Delta t_{k}} \|F_{k}\|. \label{ESTRED}
\end{align}
Therefore, from equation \eqref{ESTRED}, we obtain the estimation
\eqref{PRERED}. \qed

\vskip 2mm

In order to prove that the sequence $\{\|F(x_k)\|\}$ converges to zero when $k$
tends to infinity, we also need to estimate the lower bound of the time step
size $\Delta t_{k}$.

\vskip 2mm

\begin{lemma} \label{DTBOUND}
Assume that $F: \; \Re^{n} \to \Re^{m}$ is continuously differentiable
and its Jacobian function  $J$ is Lipschitz continuous. That is to
say, it exists a positive number $L$ such that
\begin{align}
  \|J(x)-J(y)\| \le L\|x-y\|, \; \forall x, \, y \in \Re^{n}.   \label{LIPCON}
\end{align}
Furthermore, we suppose that the sequence $\{x_k\}$ is generated by Algorithm
\ref{ALGGCNMTR} and the condition \eqref{UBINVJ} holds for all $J_{k} \,
(k = 0, \, 1, \, \ldots)$. Then, there exists a positive number
$\delta_{\Delta t}$ such that
\begin{align}
  \Delta t_{k} \ge \gamma_{2} \delta_{\Delta t} > 0  \label{DTGEPN}
\end{align}
holds for $k = 0, \, 1, \, 2, \, \ldots$, where $\Delta t_{k}$ is adaptively
adjusted by formulas \eqref{RHOK}-\eqref{TSK1}.
\end{lemma}

\vskip 2mm

\proof We prove this result by distinguishing two different cases, i.e. $J_{k}
= J(x_{k})$ or $J_{k} = J_{k-1}$. (i) Firstly, we consider the case of
$J_{k}= J(x_{k})$. From the Lipschitz continuous assumption \eqref{LIPCON} of
$J(\cdot)$, we have
\begin{align}
    & \left\|F(x_{k}+s_{k}) - F(x_k) - J(x_k)s_{k}\right\|
    = \left\|\int_{0}^{1}J(x_{k}+ ts_{k})s_{k}dt - J(x_k)s_{k}\right\|
    \nonumber \\
    & = \left\|\int_{0}^{1}(J(x_{k}+ ts_{k})-J(x_k))s_{k}dt \right\|
    \le \int_{0}^{1}\|(J(x_{k}+ ts_{k})-J(x_k))s_{k}\|dt \nonumber \\
    & \le \int_{0}^{1}\|J(x_{k}+ ts_{k})-J(x_k)\| \|s_{k}\|dt
    \le \int_{0}^{1}L \|s_{k}\|^{2} tdt = \frac{1}{2}L \|s_{k}\|^{2}.
    \label{LIPAPUB}
\end{align}
On the other hand, from equations \eqref{SMEDAE}, \eqref{UBINVJ} and
\eqref{GINVJK}, we have
\begin{align}
   \|s_k\|  = \frac{\Delta t_k}{1+ \Delta t_k} \left\|-J_{k}^{+}F_{k} \right\|
   \le \frac{\Delta t_k}{c_{\sigma}(1+ \Delta t_k)} \|F_{k}\|.  \label{ESTSK}
\end{align}
Thus, from equations \eqref{LIPAPUB}-\eqref{ESTSK}, we obtain
\begin{align}
   \left\|F(x_{k}+s_{k}) - F(x_k) - J(x_k)s_{k}\right\|
   \le \frac{L}{2c_{\sigma}^{2}}
   \left(\frac{\Delta t_k}{1+\Delta t_k}\right)^{2}\|F_{k}\|^{2}.
  \label{LINAPUPB}
\end{align}

From the definition \eqref{RHOK} of $\rho_{k}$, the estimation
\eqref{PRERED}, and equation \eqref{LINAPUPB}, we obtain
\begin{align}
   & |\rho_{k} - 1| = \left|\frac{\|F(x_k)\|-\|F(x_{k}+s_{k})\|}
   {\|F(x_k)\|-\|F(x_k)+J(x_k)s_{k}\|}-1 \right| \nonumber \\
    & \hskip 2mm \le \frac{\|F(x_{k}+s_{k}) - F(x_k) - J(x_k)s_{k}\|}
  {\|F(x_k)\|-\|F(x_k)+J(x_k)s_{k}\|}
   \le \frac{L}{2c_{\sigma}^{2}}
   \frac{\Delta t_k}{1+\Delta t_k} \|F_{k}\|
   \le \frac{L}{2c_{\sigma}^{2}} \|F_{k}\|.    \label{ESTRHO}
\end{align}
According to Algorithm \ref{ALGGCNMTR}, we know that the sequence $\{\|F(x_k)\|\}$
is monotonically decreasing. Consequently, we have $\|F(x_k)\| \le \|F(x_0)\|, \;
k = 1, \, 2, \, \dots$. We denote
\begin{align}
    \delta_{\Delta t} \triangleq \min\left\{\frac{2c_{\sigma}^{2}}{\|F(x_0)\|L} \eta_{1}, 
    \; \Delta t_{0}\right\}.   \label{PARDELAT}
\end{align}
Thus, from equations \eqref{ESTRHO}-\eqref{PARDELAT}, we obtain
$|\rho_{k} - 1| \le \eta_{1}$ when $\Delta t_{k} \le \delta_{\Delta t}$.
Consequently, according to the time-stepping scheme \eqref{TSK1},
$\Delta t_{k+1}$ will be enlarged.

\vskip 2mm

(ii) The other case is $J_{k} = J_{k-1}$. When $J_{k} = J_{k-1}$, from equation
\eqref{UPDJK1}, we know $|1-\rho_{k-1}| \le \eta_{1}$. Consequently, according
to the time-stepping scheme \eqref{TSK1}, $\Delta t_k$ will be greater than
$\Delta t_{k-1}$, i.e. $\Delta t_{k} = \gamma_{1} \Delta t_{k-1}$.

\vskip 2mm

Assume that $K$ is the first index such that $\Delta t_{K} \le \delta_{\Delta t}$.
Then, from equation \eqref{UPDJK1} and the above discussions, we know that $J_{K}
= J(x_{K})$. Otherwise, from the discussion of the case (ii), we know
$\Delta t_{K-1} < \Delta t_{K}$, which contradicts the assumption that
 $K$ is the first index such that $\Delta t_{K} \le \delta_{\Delta t}$.
 Therefore, from equations \eqref{ESTRHO}-\eqref{PARDELAT}, we obtain
$|\rho_{K} - 1| \le \eta_1$. Consequently, $\Delta t_{K+1}$ will be enlarged
according to the adaptive adjustment scheme \eqref{TSK1}. Consequently,
$\Delta t_k \ge \gamma_{2} \delta_{\Delta t}$ holds for all
$k = 0, \, 1, \, 2, \, \dots$. \qed

\vskip 2mm

By using the estimate results of Lemma \ref{LEMESTPRED} and Lemma \ref{DTBOUND}, we can prove
that the sequence $\{\|F(x_k)\|\}$ converges to zero when $k$ tends to infinity.

\begin{theorem} \label{THECOVFXK}
Assume that $F: \; \Re^{n} \to \Re^{m}$ is continuously differentiable
and its Jacobian function $J$ satisfies the Lipschitz condition \eqref{LIPCON}.
Furthermore, we suppose that the sequence $\{x_k\}$ is generated by Algorithm
\ref{ALGGCNMTR} and the Jacobian matrix $J_{k}$ satisfies the condition \eqref{UBINVJ}.
Then, we have
\begin{align}
    \lim_{k \to \infty} \inf \; \|F(x_k)\| = 0. \label{FKTOZ}
\end{align}
\end{theorem}

\proof According to Algorithm \ref{ALGGCNMTR} and Lemma \ref{DTBOUND}, we know that
there exists an infinite subsequence $\{x_{k_{l}}\}$ such that
\begin{align}
   \frac{\|F(x_{k_l})\| - \|F(x_{k_l}+s_{k_l})\|}{\|F(x_{k_l})\|-
    \|F(x_{k_l})+J_{k_l}s_{k_l}\|} \ge \eta_{2} \label{ASUBSQ}
\end{align}
holds for all $l = 0, \, 1, \, 2, \, \ldots$. Otherwise, all steps are rejected
after a given iteration index, then the time step size will keep decreasing,
which contradicts equation \eqref{DTGEPN}.

\vskip 2mm

From equations \eqref{PRERED}, \eqref{ASUBSQ} and \eqref{DTGEPN}, we have
\begin{align}
    \|F(x_{k_l})\| - \|F(x_{k_l}+s_{k_l})\| \ge \eta_{2}
    \frac{\Delta t_{k_l}}{1+\Delta t_{k_l}}\|F(x_{k_l})\|
    \ge \eta_{2} \frac{\gamma_{2}\delta_{\Delta t}}{1+ \gamma_{2}
  \delta_{\Delta t}} \|F(x_{k_l})\|. \label{FKLGE}
\end{align}
Therefore, from equation \eqref{FKLGE} and $\|F(x_{k+1})\| \le \|F(x_k)\|$, we
have
\begin{align}
    & \|F(x_0)\| \ge \|F(x_0)\| - \lim_{k \to \infty} \|F(x_k)\| =
    \sum_{k = 0}^{\infty} (\|F(x_k)\| - \|F(x_{k+1})\|) \nonumber \\
    & \hskip 2mm \ge \sum_{l = 0}^{\infty} (\|F(x_{k_l})\| - \|F(x_{k_l}+s_{k_l})\|)
    \ge \eta_{2} \frac{\gamma_{2}\delta_{\Delta t}}{1+ \gamma_{2}
  \delta_{\Delta t}} \sum_{l=0}^{\infty}\|F(x_{k_l})\|.
  \label{FX0GESUM}
\end{align}
Consequently, from equation \eqref{FX0GESUM}, we obtain
\begin{align}
    \lim_{l \to \infty} \|F(x_{k_l})\| = 0. \nonumber
\end{align}
That is to say, the result \eqref{FKTOZ} is true. \qed

\vskip 2mm

Under the full row rank of $J(x^{\ast})$ and the local Lipschitz continuity
\eqref{LIPCON}, we analyze the local superlinear convergence of Algorithm
\ref{ALGGCNMTR} near the solution $x^{\ast}$. The framework of its proof can be
roughly described as follows. Firstly, we prove that the sequence $\{x_k\}$
converges to $x^{\ast}$ when $x_{0}$ comes close enough to the solution
$x^{\ast}$. Then, we prove $\lim_{k \to \infty} \Delta t_k = \infty$. Finally,
we prove that the search step $s_{k}$ approximates the Newton step $s_{k}^{N}$.
Consequently, the sequence $\{x_{k}\}$ superlinearly converges to $x^{\ast}$.

\vskip 2mm

For convenience, we define the neighbourhood $B_{\delta}(x^{\ast})$ of $x^{\ast}$
as
\begin{align}
    B_{\delta}(x^{\ast})= \{x: \|x - x^{\ast}\|\le \delta \}. \label{NBXAST}
\end{align}

\vskip 2mm

\begin{lemma} \label{LEMCONL}
Assume that $F: \; \Re^{n} \to \Re^{m}$ is continuously differentiable and
$F(x^{\ast}) = 0$. Furthermore, we suppose that its Jacobian function $J$
satisfies the Lipschitz continuity \eqref{LIPCON} and the condition \eqref{UBINVJ}
when $x \in B_{\delta}(x^{\ast})$. Then, there exists a neighborhood
$B_{r}(x^{\ast})$ of $x^{\ast}$ such that the sequence $\{x_k\}$ generated by
Algorithm \ref{ALGGCNMTR} with $x_{0} \in B_{r}(x^{\ast})$ converges
to $x^{\ast}$.
\end{lemma}
\proof From equations \eqref{SVD}-\eqref{UBINVJ}, we obtain the generalized
inverse $J_{k}^{+}$ in equation \eqref{GINVJK} and its estimation
\begin{align}
   \left\|J_{k}^{+}\right\| \le {1}/{c_{\sigma}}, \; k = 0, \, 1, \, 2, \ldots.
   \label{INJACUPB}
\end{align}
We denote $e_{k} = x_{k} - x^{\ast}$. When $s_{k}$ is not an accepted step, we
obviously have $e_{k+1} = e_{k}$. Therefore, we consider the case that $s_{k}$
is an accepted step. When $s_{k}$ is an accepted step, from the generalized
continuation Newton method \eqref{GICNM}, we have
\begin{align}
   e_{k+1} & = e_{k} + s_{k} = e_{k} + \frac{\Delta t_{k}}{1+ \Delta t_{k}}
   J_{k}^{+}(F(x_k) - F(x^{\ast}))
   \nonumber \\
   & = e_{k} + \frac{\Delta t_{k}}{1+ \Delta t_{k}}
   J_{k}^{+}\int_{0}^{1} J(x^{\ast} + t e_{k})e_{k}dt.    \label{EKSK}
\end{align}
By rearranging the above equation \eqref{EKSK}, we obtain
\begin{align}
   e_{k+1} = \frac{1}{1+\Delta t_k}e_{k} + \frac{\Delta t_{k}}{1+ \Delta t_{k}}
    J_{k}^{+}\int_{0}^{1} \left(J(x^{\ast} + t e_{k}) - J(x_k) \right)e_{k}dt.
    \nonumber
\end{align}
By using the Lipschitz continuity \eqref{LIPCON} of $J$ and the
estimation \eqref{INJACUPB}, we have
\begin{align}
   & \|e_{k+1}\|  \le \frac{1}{1+\Delta t_k}\|e_k\|
   + \frac{\Delta t_k}{1+ \Delta t_k} \left\| J_{k}^{+} \right\| \int_{0}^{1}
     \left\|J(x^{\ast} + t e_{k}) -J(x_k)\right\| \|e_{k}\|dt
     \nonumber \\
   & \hskip 2mm  \le \frac{1}{1+\Delta t_k}\|e_k\| + \frac{\Delta t_k}{1+ \Delta t_k}
   \frac{L}{2c_{\sigma}}\|e_k\|^{2}
   = \frac{1+ L/(2c_{\sigma})\|e_{k}\|{\Delta t_k}} {1+\Delta t_k}\|e_k\|.
      \label{EK1EKNS}
\end{align}

\vskip 2mm

We denote
\begin{align}
   q_{k} \triangleq \frac{1+ L/(2c_{\sigma})\|e_{k}\|{\Delta t_k}}
   {1+\Delta t_k}, \label{QK}
\end{align}
and select $x_{0} \in B_{\delta}(x^{\ast})$  to satisfy
\begin{align}
  \|e_{0}\| < \frac{c_{\sigma}}{L}.   \label{EKLESS1}
\end{align}
We denote $r = \min \{\delta, \, 2c_{\sigma}/L\}$.  When $x_{0} \in B_{r}(x^{\ast})$,
from equations \eqref{EK1EKNS}-\eqref{EKLESS1}, by induction, we have
\begin{align}
   \|e_{k+1}\| \le q_{k} \|e_k\|, \;
    q_{k} < \frac{1+ (1/2) \Delta t_k}{1+\Delta t_k} < 1.
    \label{CONEK1}
\end{align}

\vskip 2mm

It is not difficult to know that $f(t) \triangleq (1+\alpha t)/(1+t)$ is monotonically
decreasing when $0 \le \alpha < 1$. Thus, from the estimation \eqref{DTGEPN} of
the time step size $\Delta t_{k}$ and equation \eqref{CONEK1}, we obtain
\begin{align}
   \|e_{k+1}\| \le q_{k} \|e_{k}\| \le q \|e_{k}\|, \;
   q \triangleq \frac{1+ (1/2) \gamma_{2} \delta_{\Delta t}}
   {1+\gamma_{2}\delta_{\Delta t}} < 1.     \label{AEKLEQEK1}
\end{align}

\vskip 2mm

Consequently, from equation \eqref{AEKLEQEK1}, we know that
$e_{k+1} \le e_{k}$ holds for all $k = 0, \, 1, \, 2, \, \ldots$, since
$e_{k+1} = e_{k}$ when $s_{k}$ is not an accepted step. According to Algorithm
\ref{GCNMTR} and Lemma \ref{DTBOUND}, we know that there exists an infinite
subsequence $\{x_{k_{l}}\}$ such that $s_{k_l} \, (l = 0, \, 1, \, \ldots)$ are
all accepted steps. Otherwise, all steps are rejected after a given iteration
index, then the time step size will keep decreasing, which contradicts equation
\eqref{DTGEPN}. Therefore, from equation \eqref{AEKLEQEK1} and $e_{k} \le e_{k+1}$,
we have
\begin{align}
   e_{k_l} \le q \, e_{k_{l-1}} \le \cdots \le q^{l} e_{k_0}. \nonumber
\end{align}
That is to say, we have $\lim_{l \to \infty} e_{k_l} = 0$. By combining it with
$e_{k+1} \le e_{k}$, we obtain $\lim_{k \to \infty} x_k = x^{\ast}$. \qed

\vskip 2mm

\begin{lemma} \label{LEMTSTI}
Assume that $F: \; \Re^{n} \to \Re^{m}$ is continuously differentiable and
$F(x^{\ast}) = 0$. Furthermore, we suppose that its Jacobian function $J$
satisfies the Lipschitz continuity \eqref{LIPCON} and the condition \eqref{UBINVJ}
when $x \in B_{\delta}(x^{\ast})$. Then, there exists a neighborhood
$B_{r}(x^{\ast})$ of $x^{\ast}$ such that the sequence $\{x_k\}$ generated by
Algorithm \ref{ALGGCNMTR} with $x_{0} \in B_{r}(x^{\ast})$ converges
to $x^{\ast}$ and the generated time step $\Delta t_{k}$ tends to infinity.
\end{lemma}
\proof The first part of the lemma is proved in Lemma \ref{LEMCONL}, i.e.
$\lim_{k \to \infty} x_{k} = x^{\ast}$. Now, we prove the second part of the
lemma, i.e. $\lim_{k \to \infty} \Delta t_{k} = \infty$.

\vskip 2mm

We can assume that there exists an infinite subsequence  $\{x_{k_l}\}$
such that $J_{k_l} = J(x_{k_l})$ holds for all $l = 0, \, 1, \ldots$.
Otherwise, according to equation \eqref{UPDJK1}, all Jacobian matrices
$J_{k} \, (k = K+1, \, K+2, \, \ldots)$ equal $J_{K}$ and
$|\rho_{k} - 1| \le \eta_{1} \, (k = K, \, K+2, \, \ldots)$ after a
given iteration index $K$. Then, according to the time-stepping scheme
\eqref{TSK1}, we obtain $\Delta t_{k+1} = \gamma_{1} \Delta t_{k} \,
(k = K, \, K+1, \, \ldots)$. Consequently, we have
$\lim_{k \to \infty} \Delta t_{k} = \infty$. That is to say, for this
case, the second part of the lemma also is proved.

\vskip  2mm

Since $J_{k_l} = J(x_{k_l})$, from  equations \eqref{GICNM} and \eqref{INJACUPB},
we have
\begin{align}
  & \|s_{k_l}\|  =  \frac{\Delta t_{k_l}}{1+\Delta t_{k_l}}
  \left\|J_{k_l}^{+}F(x_{k_l})\right\|
  \le \frac{\Delta t_{k_l}}{1+\Delta t_{k_l}} \|J_{k_l}^{+}\| \|F(x_{k_l})\| \nonumber \\
  & \hskip 2mm \le \frac{\Delta t_{k_l}}{c_{\sigma}(1+\Delta t_{k_l})}\|F(x_{k_l})\|. \label{SKUPB}
\end{align}
Similarly to the estimation \eqref{ESTRHO}, from the
definition \eqref{RHOK} of $\rho_{k_l}$, inequalities \eqref{PRERED} and
\eqref{SKUPB}, we have
\begin{align}
  & |\rho_{k_l} - 1|  = \left|\frac{\|F(x_{k_l})\|-\|F(x_{k_l}+s_{k_l})\|}
  {\|F(x_{k_l})\|-\|F(x_{k_l})+J(x_{k_l})s_{k_l}\|}-1 \right| \nonumber \\
   & \hskip 2mm \le \frac{L}{2c_{\sigma}^{2}}
   \left(\frac{\Delta t_{k_l}}{1+\Delta t_{k_l}}\right) \|F(x_{k_l})\|
   \le \frac{L}{2c_{\sigma}^{2}} \|F(x_{k_l})\|.  \label{ESTRHOEN}
\end{align}
Since $\lim_{k \to \infty} x_k =  x^{\ast}$ and $F(x^{\ast}) = 0$, we can
select a sufficiently large number $K$ such that
\begin{align}
   \|F(x_{K})\| \le \frac{2\eta_{1}c_{\sigma}^{2}}{3L}.  \label{FXKLES}
\end{align}
From inequalities \eqref{ESTRHOEN}-\eqref{FXKLES} and the monotonically
decreasing property $\|F(x_{k+1})\| \le \|F(x_{k})\|$, we have
$|\rho_{k_l} - 1| \le \eta_1$ when $k_l \ge K$. This means $\Delta t_{{k_l}+1}
= \gamma_{1} \Delta t_{k_l}$ according to the time-stepping scheme \eqref{TSK1}.

\vskip 2mm

Now, we consider the $(k_{l}+1)$-th iteration. From equation \eqref{UPDJK1}, we know
that $J_{{k_l}+1} = J_{k_l} = J(x_{k_l})$. Then, from the definition \eqref{RHOK} of
$\rho_{{k_l}+1}$, equation \eqref{PRERED} and the Lipschitz continuity
\eqref{LIPCON}, we have
\begin{align}
   & |\rho_{{k_l}+1} - 1|  = \left|\frac{\|F(x_{{k_l}+1})\|-\|F(x_{{k_l}+1}+s_{{k_l}+1})\|}
    {\|F(x_{{k_l}+1})\|-\|F(x_{{k_l}+1})+J_{{k_l}+1}s_{{k_l}+1}\|}-1 \right| \nonumber \\
    & \hskip 2mm \le \frac{\|F(x_{{k_l}+1}+s_{{k_l}+1})
     - F(x_{{k_l}+1})  - J_{{k_l}+1}s_{{k_l}+1}\|}
    {\|F(x_{{k_l}+1})\|-\|F(x_{{k_l}+1})+J_{{k_l}+1}s_{{k_l}+1}\|} \nonumber \\
    & \hskip 2mm = \frac{\|\int_{0}^{1} (J(x_{{k_l}+1}
    + ts_{{k_l}+1}) - J(x_{k_l}))s_{{k_l}+1}dt\|}
    {\|F(x_{{k_l}+1})\|-\|F(x_{{k_l}+1})+J_{{k_l}+1}s_{{k_l}+1}\|} \nonumber \\
   & \hskip 2mm \le \frac{1+\Delta t_{{k_l}+1}}{\Delta t_{{k_l}+1}}
   \frac{L \left(0.5 \|s_{{k_l}+1}\|^{2}+\|s_{{k_l}+1}\| \|s_{k_l}\|\right)}
   {\|F(x_{{k_l}+1})\|}.      \label{UNESTRHOEN}
\end{align}
By substituting equation \eqref{SKUPB} into equation \eqref{UNESTRHOEN}, we obtain
\begin{align}
     |\rho_{{k_l}+1} - 1| \le \frac{L \; \Delta t_{{k_l}+1}}
     {c_{\sigma}^{2}(1+\Delta t_{{k_l}+1})}
     \left(0.5\|F(x_{{k_l}+1})\|+\|F(x_{k_l})\|\right)
     \le \frac{3L}{2c_{\sigma}^{2}} \|F(x_{k_l})\|, \label{ROHKNAP}
\end{align}
where the property $\|F(x_{{k_l}+1})\| \le \|F(x_{k_l})\|$ is used in the last
inequality.

\vskip 2mm

From equations \eqref{FXKLES} and \eqref{ROHKNAP}, we have
$|\rho_{{k_l}+1} - 1| \le \eta_1$. This means $\Delta t_{{k_l}+2}
= \gamma_{1} \Delta t_{{k_l}+1} = \gamma_{1}^{2} \Delta t_{k_l}$ according to
the time-stepping scheme \eqref{TSK1}. Thus, for the $(k_{l}+2)$ iteration, when
$|\rho_{{k_l}+2} - 1| > \eta_{1}$, according to the time-stepping scheme
\eqref{TSK1}, we have $\Delta t_{{k_l}+3} \ge \gamma_{2} \gamma_{1}^{2}
\Delta t_{k_l} = \gamma_{1} \Delta t_{k_l}$ (Here, we select $\gamma_{1} = 2$
and $\gamma_{2} = 1/2$). Furthermore, from equation \eqref{UPDJK1}, we know
$J_{{k_l}+3} = J(x_{{k_l}+3})$ at the $({k_l}+3)$-th iteration. Similarly to
the estimation $|\rho_{k_l} - 1| \le \eta_{1}$ of $J_{k_l} = J(x_{k_l})$, we
have $|\rho_{{k_l}+3} - 1| \le \eta_1$. This means $\Delta t_{{k_l}+4}
= \gamma_{1} \Delta t_{{k_l}+3} \ge \gamma_{1}^{2} \Delta t_{k_l}$ according
to the time-stepping scheme \eqref{TSK1}. Thus, the subsequent iterations start
a new cycle for the time step size.

\vskip 2mm

When $|\rho_{{k_l}+2} - 1| \le \eta_{1}$, according to the
time-stepping scheme \eqref{TSK1}, we have $\Delta t_{{k_l}+3} = \gamma_{1}
\Delta t_{{k_l}+2} = \gamma_{1}^{3} \Delta t_{k_l}$ and the time steps keep
increasing until $|\rho_{{k_l}+m} - 1| > \eta_{1} \, (m = 3, \, 4, \, \ldots)$.
Then, the subsequent iterations start a new cycle for the time step size.

\vskip 2mm

By combining the above discussions of two cases $|\rho_{{k_l}+2} - 1| > \eta_{1}$
or $|\rho_{{k_l}+2} - 1| \le \eta_{1}$, we know $\Delta t_{k_l} \ge
\gamma_{1} \Delta t_{k_{l+1}}$ and $\Delta t_{{k_l}+m} \ge \Delta t_{k_l}$
when $ 0 \le m \le k_{l+1} - k_{l}$. Consequently, we obtain $\lim_{l \to \infty}
\Delta t_{k_l} = \infty$. By combining the property $\Delta t_{{k_l}+m} \ge
\Delta t_{k_l} \, (0 \le m \le k_{l+1} - k_{l})$, we obtain
$\lim_{k \to \infty} \Delta t_{k} = \infty$. \qed

\vskip 2mm

\begin{theorem} \label{LOCSCON}
Assume that $F: \; \Re^{n} \to \Re^{m}$ is continuously differentiable and
$F(x^{\ast}) = 0$. Furthermore, we suppose that its Jacobian function $J$
satisfies the Lipschitz continuity \eqref{LIPCON} and the condition \eqref{UBINVJ}
when $x \in B_{\delta}(x^{\ast})$. Then, there exists a neighborhood
$B_{r}(x^{\ast})$ of $x^{\ast}$ such that the sequence $\{x_k\}$ generated by
Algorithm \ref{ALGGCNMTR} with $x_{0} \in B_{r}(x^{\ast})$ converges
superlinearly to $x^{\ast}$.
\end{theorem}
\proof From Lemma \ref{LEMCONL} and Lemma \ref{LEMTSTI}, we know
$\lim_{k \to \infty} x_{k}  = x^{\ast}$ and $\lim_{k \to \infty} \Delta t_{k}
= \infty$. Firstly, we prove that there are only finite steps which are rejected.
That is to say, all steps are accepted after a given iteration index.

\vskip 2mm

We assume that there exist the infinite rejected steps. Since
$\lim_{k \to \infty} x_k =  x^{\ast}$ and $F(x^{\ast}) = 0$, we can
select a sufficiently large number $K_{1}$ such that
\begin{align}
   \|F(x_{K_1})\| \le \frac{\eta_{1}c_{\sigma}^{2}}{L}.  \label{FXKLENA}
\end{align}
Furthermore, there exists a positive $K_{2}$ such that
\begin{align}
      \|x_{k} - x_{l}\| \le \frac{\eta_{1} c_{\sigma}}{2L},
      \; \forall k, \, l \ge K_{2}.  \label{CAUCHYSQ}
\end{align}
We denote $K = \max\{K_{1}, \, K_{2}\}$. For the $k$-th iteration, we assume
that $s_{l-1} \, (l \le k) $ is the step such that $|\rho_{l-1} - 1|
< \eta_{1}$ holds and its index $l$ is the closest to $k$. Then, we have
$J_{k} = J_{l} = J(x_{l})$ according to equation \eqref{UPDJK1}.

\vskip 2mm

Similarly to the estimation \eqref{UNESTRHOEN}, from the definition \eqref{RHOK}
of $\rho_{k}$, equation \eqref{GICNM} and the Lipschitz continuity \eqref{LIPCON},
we have
\begin{align}
     |\rho_{k} - 1| \le \frac{1+\Delta t_{k}}{\Delta t_{k}}
   \frac{L \left(0.5 \|s_{k}\|^{2}+\|s_{k}\| \|x_{k} - x_{l}\|\right)}
   {\|F(x_{k})\|}. \label{RHOUALK}
\end{align}
By substituting equation \eqref{ESTSK} into equation \eqref{RHOUALK}, we obtain
\begin{align}
     & |\rho_{k} - 1| \le \frac{L}{c_{\sigma}}\left(\frac{1}{2c_{\sigma}}
     \frac{\Delta t_{k}}{1+\Delta t_{k}}\|F(x_{k})\| + \|x_{k} - x_{l}\|\right)
     \nonumber \\
     & \hskip 2mm \le \frac{L}{c_{\sigma}}\left(\frac{1}{2c_{\sigma}}
     \|F(x_{k})\| + \|x_{k} - x_{l}\|\right). \label{RHOKNALE}
\end{align}

\vskip 2mm

By substituting inequalities \eqref{FXKLENA}-\eqref{CAUCHYSQ} and the
monotonically decreasing property $\|F(x_{k+1})\| \le \|F(x_{k})\|$ into
equation \eqref{RHOKNALE}, we have $|\rho_{k} - 1| \le \eta_1$ when
$k \ge K$. This means $\Delta t_{k+1} = \gamma_{1} \Delta t_{k}$ according
to the time-stepping scheme \eqref{TSK1}. Thus, we know that all steps
$s_{k} \, (k \ge K)$ are the accepted steps, which contradicts the assumption
of the infinite rejected steps. Therefore, there exist only finite rejected steps.

\vskip 2mm

We denote $e_{k} = x_{k} - x^{\ast}$. Then, similarly to the estimation
\eqref{EK1EKNS}, from the Lipschitz continuity \eqref{LIPCON} and the
estimation \eqref{INJACUPB}, we have
\begin{align}
     & \frac{\|e_{k+1}\|}{\|e_{k}\|}   \le
     \frac{1}{1+\Delta t_k}  + \frac{\Delta t_k}{1+ \Delta t_k}
     \frac{L}{2c_{\sigma}}\|e_k\|
     \le \frac{1}{1+\Delta t_k} + \frac{L}{2c_{\sigma}}\|e_k\|, \; k \ge K.
     \label{EK1VSEK}
\end{align}
By substituting $\lim_{k \to \infty} \Delta t_k = \infty$ and
$\lim_{k \to \infty} \|e_{k}\| = 0$ into equation \eqref{EK1VSEK}, we
obtain
\begin{align}
   \lim_{k \to \infty} \frac{\|e_{k+1}\|}{\|e_{k}\|} = 0. \nonumber
\end{align}
That is to say, the sequence $\{x_{k}\}$ superlinearly converges to $x^{\ast}$.
\qed

\vskip 2mm

\section{Numerical experiments}

\vskip 2mm

Since the classical Homotopy continuation method such as HOMPACK90
\cite{WSMMW1997} can not effectively tackle the underdetermined system of
nonlinear equations, we only compare Algorithm \ref{ALGGCNMTR} (GCNMTr) with the
traditional optimization method such as the Levenberg-Marquardt method
(the built-in subroutine fsolve.m of the MATLAB R2020a environment
\cite{FY2005,ISU2019,Levenberg1944,LLT2007,MATLAB,Marquardt1963,More1978}).
The Jacobian matrix $J(x)$ of Algorithm \ref{ALGGCNMTR} is approximated by 
the difference formula \eqref{NUMHESS}. The codes are executed by a HP 
notebook with the Intel quad-core CPU and 8Gb memory in the MATLAB R2020a 
environment \cite{MATLAB}. 

\vskip 2mm

At every iteration, in order to obtain the accepted trial step $s_{k}$, the
Levenberg-Marquardt method \cite{MATLAB,More1978} needs to solve several linear
least-squares problems as follows:
\begin{align}
     \begin{bmatrix} J_{k} \\ \sqrt{\lambda}D_{k} \end{bmatrix}
      \approx - \begin{bmatrix} F_{k} \\ 0 \end{bmatrix}, \label{AUGLLS}
\end{align}
such that $\|s_{k}(\lambda)\| \approx \Delta_{k}$, where $\Delta_{k}$ is
the trust-region radius and the scaled matrix $D_{k}$  is usually selected
as a diagonal matrix as follows \cite{More1978}:
\begin{align}
     & D_{k} = \text{diag}\left(d_{1}^{(k)}, \, \ldots, \, d_{n}^{(k)}\right), \nonumber \\
     & d_{i}^{(k)} = \max\{d_{i}^{(k-1)}, \, \|J(x_{k})_{i}\|\},  \,
     i = 1, \, 2, \, \ldots, \, n. \nonumber
\end{align}

\vskip 2mm

The test underdetermined problems of nonlinear equations are derived from
\cite{Andrei2008,Luksan1994,MGH1981,SB2013}. We preserve the first $m$ elements of the
gradient $g(x)$ of the unconstrained optimization function $f(\cdot)$ as the test
underdetermined system, i.e.
\begin{align}
   F(x) = [g_{1}(x), \, \ldots, \, g_{m}(x)]^{T}, \nonumber
\end{align}
where $g(x) = \nabla f(x), \, x \in \Re^{n}$. Their dimensions are all set by
$(n = 2000, \;  m = n-1)$, $(n = 2000, \;  m = 10)$ and $(n = 2000, \;  m = n)$. 
The initial point of the test problem is set by
$x_{0} = \text{ones} (n, \,  1)$ when $\text{ones} (n, \,  1)$ is not the zero point
of $F(x)$. Otherwise, the initial point is set  $x_{0} = 2 \times \text{ones} (n, \, 1)$.
The tolerable errors of two methods are both set by
\begin{align}
   \|F(x^{it})\|_{\infty} \le 10^{-6}. \label{TERMCON}
\end{align}

\vskip 2mm

\begin{figure}[h]
   \centering
   \includegraphics[width=0.80\textwidth]{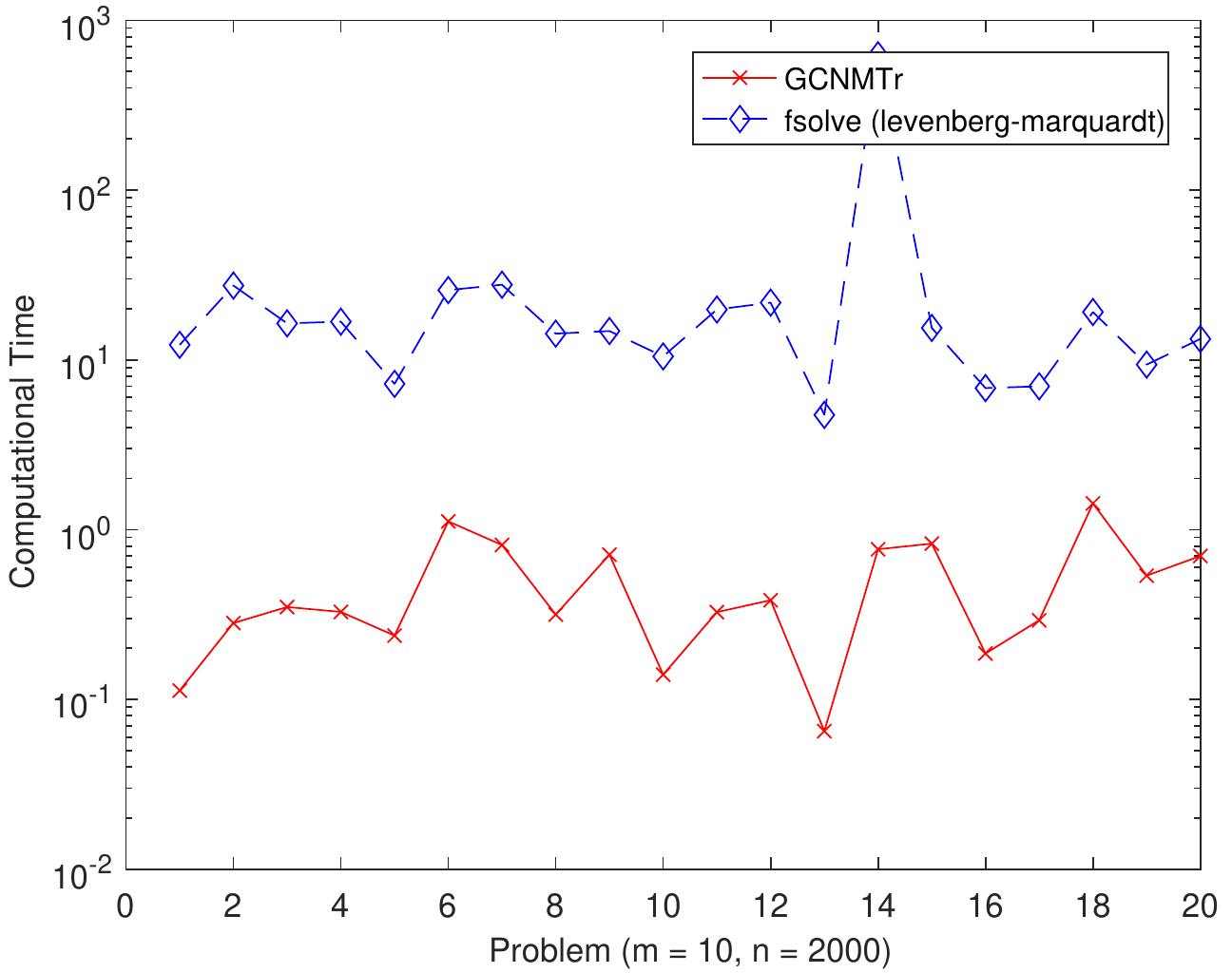}
   \caption{The computational time of GCNMTr and fsolve for the underdetermined problems with 
   $m = 10, \,  n = 2000$.}
   \label{FIGTIME1}
\end{figure}

\begin{figure}[h]
   \centering
   \includegraphics[width=0.80\textwidth]{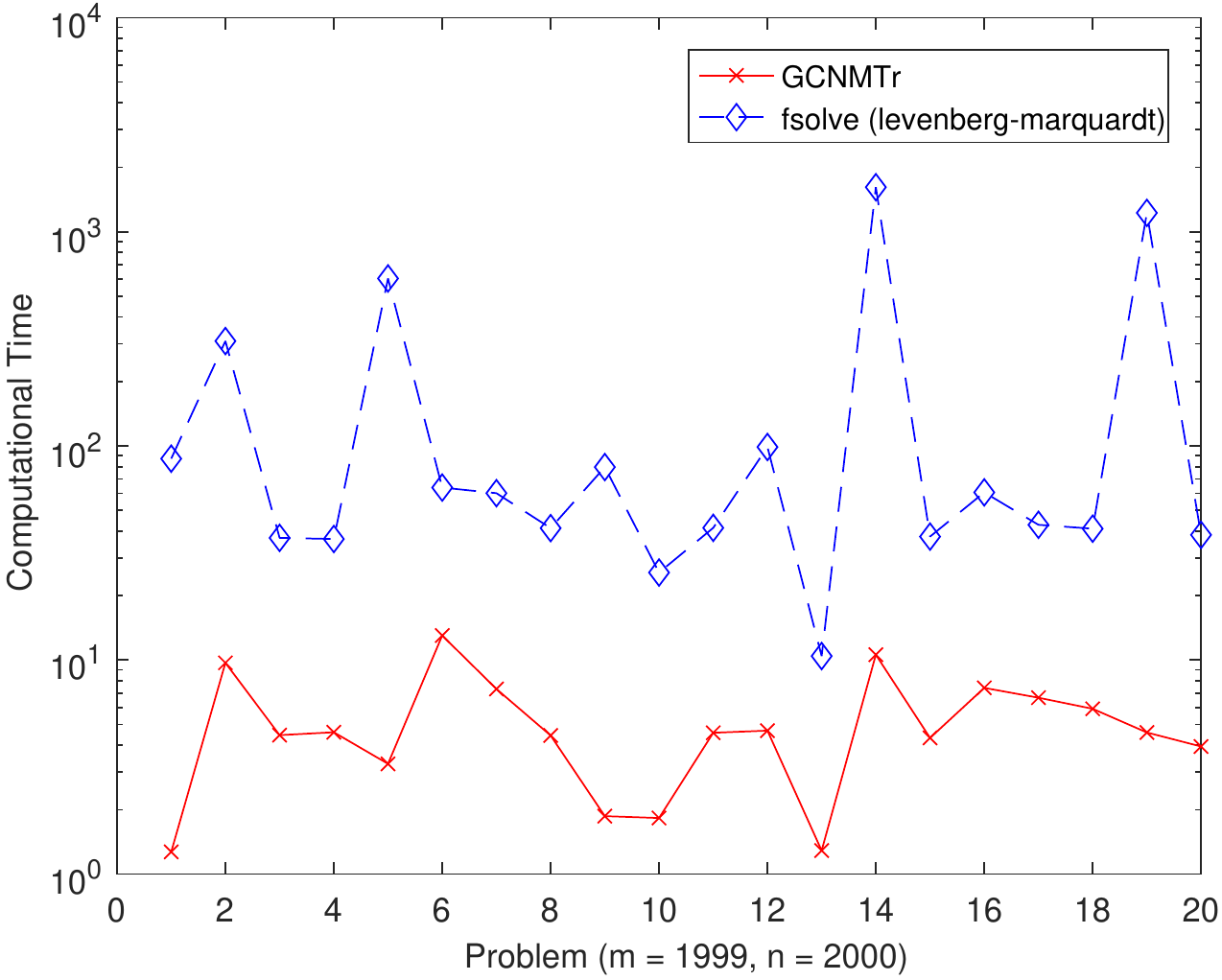}
   \caption{The computational time of GCNMTr and fsolve for the underdetermined problems
   with $m = 1999, \;  n = 2000$.}
   \label{FIGTIME2}
\end{figure}

\begin{figure}[h]
   \centering
   \includegraphics[width=0.80\textwidth]{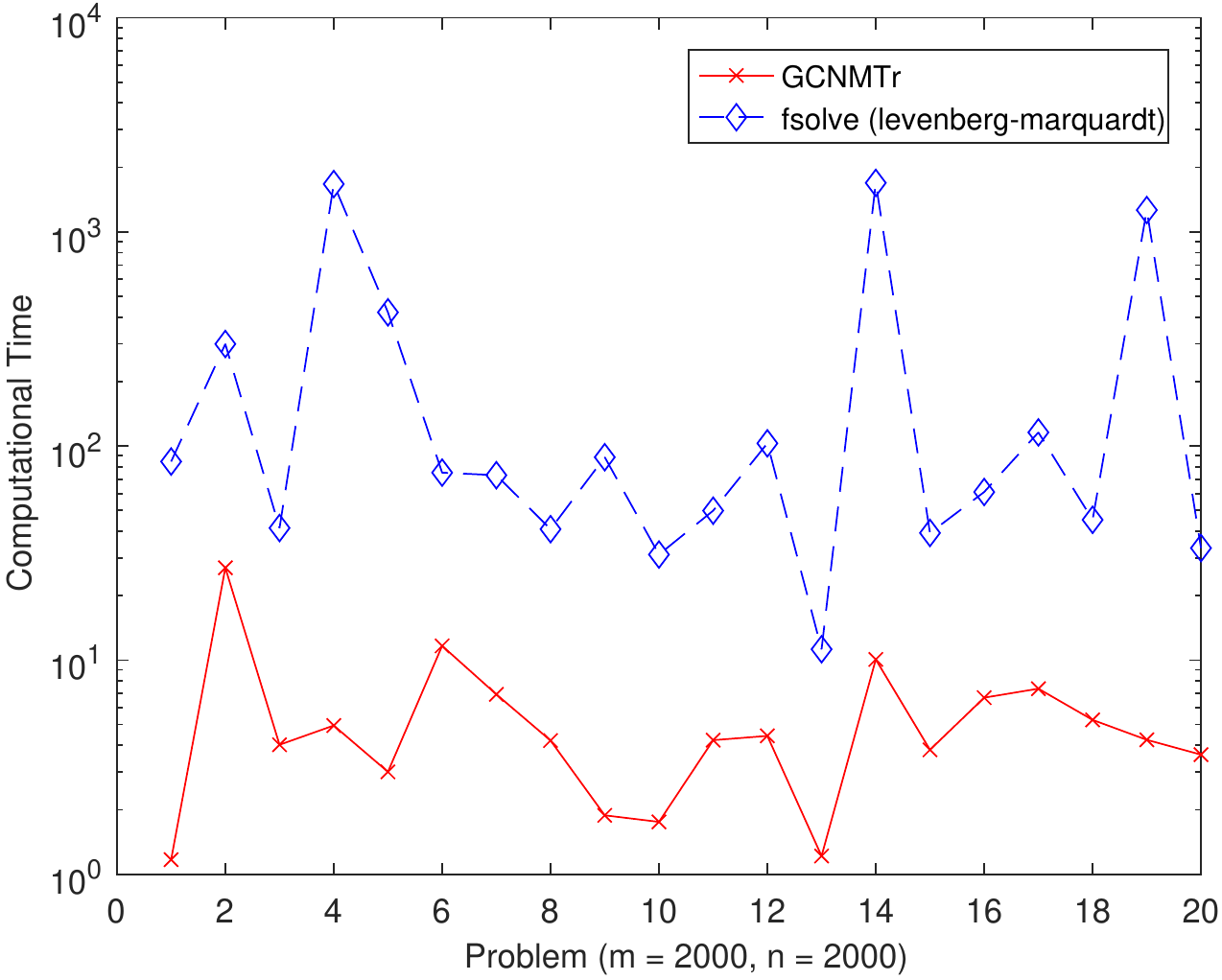}
   \caption{The computational time of GCNMTr and fsolve for the determined problems with 
   $m = 2000, \; n = 2000$.}
   \label{FIGTIME3}
\end{figure}

The numerical results are arranged in Tables \ref{TABCOMT1}-\ref{TABCOMT3}.
NJ stands for the number of the Jacobian evaluations required for convergence
in Tables \ref{TABCOMT1}-\ref{TABCOMT3}. The computational time  of GCNMTr and 
fsolve is illustrated by Figures \ref{FIGTIME1}-\ref{FIGTIME3}. From Tables 
\ref{TABCOMT1}-\ref{TABCOMT3}, we find that GCNMTr performs well for those test 
problems, and the levenberg-marquardt method (fsolve) fails to solve some problems.

\vskip 2mm

Furthermore, from Figures \ref{FIGTIME1}-\ref{FIGTIME3} and Tables 
\ref{TABCOMT1}-\ref{TABCOMT3}, we also find that GCNMTr is faster than the 
Levenberg-Marquardt method (fsolve) and the computational time of GCNMTr is 
about $1/8$ to $1/50$ of that of fsolve.  One of the reasons is that GCNMTr 
does not need to update the Jacobian matrix $J_{k}$ and decompose it when it 
performs well. This strategy can save much computational time, in comparison 
to that of the updating the Jacobian matrix $J(x_{k})$ at every iteration for 
the traditional Levenberg-Marquardt method. The other
reason is that fsolve uses the QR decomposition to solve the linear system
\eqref{AUGLLS}, which requires $2n^{2}(2n/3+m)$ flops (p. 264, \cite{GV2013}).
However, GCNMTr uses the QR decomposition to solve the linear systems of equation
\eqref{SKJTQR}, which only requires  $2m^{2}(n - m/3)$ flops and
about the half of that of fsolve.

\vskip 2mm

\begin{table}[!http]
  \newcommand{\tabincell}[2]{\begin{tabular}{@{}#1@{}}#2\end{tabular}}
  \scriptsize
  \centering
  \caption{Numerical results of GCNMTr and fsolve for the underdetermined problems with $m = 10, \;  n = 2000$.}
  \label{TABCOMT1}
  \resizebox{\textwidth}{!}{
  \begin{tabular}{|c|c|c|c|c|c|c|c|c|c|}
  \hline
  \multirow{2}{*}{Problems } & \multicolumn{2}{c|}{GCNMTr} & \multicolumn{2}{c|}{fsolve.m (levenberg-marquardt)}  \\ \cline{2-5}
  & \tabincell{c}{NJ (Iteration, time$/s$)}     & $||F(x^{it})||_\infty$  & \tabincell{c}{Iteration (time$/s$)} & $||F(x^{it})||_\infty$     \\\hline

  \tabincell{c}{1. Trid Function \cite{SB2013} \\ (m = 10, n = 2000)}   & \tabincell{c}{2 (14, 0.1127$s$)} & 1.6495E-08 & \tabincell{c}{5 (12.2854$s$)} & 2.8635E-11   \\ \hline

  \tabincell{c}{2. Grewank Function \cite{SB2013} \\ (m = 10, n=2000)}   & \tabincell{c}{2 (11, 0.2816$s$)} & 5.6776E-07 & \tabincell{c}{10 (27.3960$s$)} & 5.0238E-16 \\ \hline

  \tabincell{c}{3. Dixon Price Function \cite{SB2013} \\ (m = 10, n = 2000)}   & \tabincell{c}{5 (20, 0.3500$s$)} & 5.2291E-07 & \tabincell{c}{7 (16.4054$s$)} & 1.2434E-14 \\ \hline

  \tabincell{c}{4. Rosenbrock Function \cite{MGH1981} \\ (m = 10, n = 2000)}   & \tabincell{c}{5 (24, 0.3274$s$)} & 4.1377E-07 & \tabincell{c}{7 (16.7926$s$)} & 6.1729E-11  \\ \hline

  \tabincell{c}{5. Trigonometric Function  \cite{MGH1981} \\ (m = 10, n = 2000)}   & \tabincell{c}{2 (11, 0.2375$s$)} & 1.4611E-07 & \tabincell{c}{3 (7.2162$s$)} & 8.5720E-12  \\ \hline

  \tabincell{c}{6. Singular Broyden Function \cite{Luksan1994} \\ (m = 10, n = 2000)}   & \tabincell{c}{19 (35, 1.1180$s$)} & 4.3241E-07  & \tabincell{c} {\red {12 (25.7608$s$)}}  & \tabincell{c} { \red{0.0213} \\ \red{(failed)}}\\ \hline

  \tabincell{c}{7. Extended Powell Singular Function \cite{MGH1981} \\ (m = 10, n = 2000)}   & \tabincell{c}{10 (25, 0.8116$s$)} & 7.7570E-07 & \tabincell{c}{12 (27.7080$s$)} & 7.5422E-07  \\ \hline

  \tabincell{c}{8. Tridiagonal System Function \cite{Luksan1994} \\ (m = 10, n = 2000)}   & \tabincell{c}{6 (20, 0.3157$s$)} & 4.8995E-07 & \tabincell{c}{8 (14.2731$s$)} & 2.2095E-08  \\ \hline

  \tabincell{c}{9. Discrete Boundary-Value Function \cite{Luksan1994} \\ (m = 10, n = 2000)}   & \tabincell{c}{2 (13, 0.7122$s$)} & 6.8392E-07 & \tabincell{c}{5 (14.7863$s$)} & 7.5662E-13 \\ \hline

  \tabincell{c}{10. Broyden Tridiagonal Function \cite{Luksan1994} \\ (m = 10, n = 2000)}   & \tabincell{c}{2 (15, 0.1398$s$)} & 1.9205E-07 & \tabincell{c}{5 (10.4810$s$)} & 1.5421E-10 \\ \hline

  \tabincell{c}{11. Extended Wood Function \cite{Andrei2008} \\ (m = 10, n = 2000)}   & \tabincell{c}{6 (22, 0.3264$s$)} & 7.0474E-08 & \tabincell{c}{8 (19.7991$s$)} & 5.8124E-09  \\ \hline

  \tabincell{c}{12. Extended Cliff Function \cite{Andrei2008} \\ (m = 10, n = 2000)}   & \tabincell{c}{6 (21, 0.3839$s$)} & 9.6747E-07 & \tabincell{c}{10 (21.7415$s$)} & 1.5466E-08 \\ \hline

  \tabincell{c}{13. Extended Hiebert Function \cite{Andrei2008} \\ (m = 10, n = 2000)}   & \tabincell{c}{2 (13, 0.0650$s$)} & 1.4634E-08 & \tabincell{c}{2 (4.7399$s$)} & 5.0164E-12  \\ \hline

  \tabincell{c}{14. Extended Maratos Function \cite{Andrei2008} \\ (m = 10, n = 2000)}   & \tabincell{c}{16 (37, 0.7660$s$)} & 7.4567E-07 & \tabincell{c} {\red {199 (623.0145$s$)}  }& \tabincell{c} {\red{0.1757} \\ \red{(failed)}} \\ \hline

  \tabincell{c}{15. Extended Psc1 Function \cite{Andrei2008} \\ (m = 10, n = 2000)}   & \tabincell{c}{5 (20, 0.8269$s$)} & 6.1415E-07 & \tabincell{c}{7 (15.4128$s$)} & 2.3426E-14  \\ \hline

  \tabincell{c}{16. Extended Quadratic Penalty QP 1 \\ Function \cite{Andrei2008}  (m = 10, n = 2000)}   & \tabincell{c}{2 (18, 0.1863$s$)} & 6.1455E-08 & \tabincell{c}{3 (6.8182$s$)} & 1.1303E-14 \\ \hline

  \tabincell{c}{17. Extended Quadratic Penalty QP 2  \\ Function \cite{Andrei2008} (m = 10, n = 2000)}   & \tabincell{c}{2 (20, 0.2924$s$)} & 5.7408E-07 & \tabincell{c}{3 (6.9848$s$)} & 6.5782E-12 \\ \hline

  \tabincell{c}{18. Extended TET Function \cite{Andrei2008} \\ (m = 10, n = 2000)}   & \tabincell{c}{7 (22, 1.4268$s$)} & 2.1258E-07 & \tabincell{c}{8 (19.1292$s$)} & 7.8478E-11 \\ \hline

  \tabincell{c}{19. EG2 Function \cite{Andrei2008} \\ (m = 10, n = 2000)}   & \tabincell{c}{2 (16, 0.5359$s$)} & 6.3057E-08 & \tabincell{c}{4 (9.3791$s$)} & 6.1195E-13 \\ \hline

  \tabincell{c}{20. Extended BD1 Function \cite{Andrei2008} \\ (m = 10, n = 2000)}   & \tabincell{c}{5 (19, 0.6993$s$)} & 1.1822E-07 & \tabincell{c}{6 (13.2805$s$)} & 1.4782E-08 \\ \hline

\end{tabular}}
\end{table}

\vskip 2mm

\begin{table}[!http]
  \newcommand{\tabincell}[2]{\begin{tabular}{@{}#1@{}}#2\end{tabular}}
  \scriptsize
  \centering
  \caption{Numerical results of GCNMTr and fsolve for the underdetermined problems 
  with $m = 1999, \;  n = 2000$.}
  \label{TABCOMT2}
  \resizebox{\textwidth}{!}{
  \begin{tabular}{|c|c|c|c|c|c|c|c|c|c|}
  \hline
  \multirow{2}{*}{Problems } & \multicolumn{2}{c|}{GCNMTr} & \multicolumn{2}{c|}{fsolve.m (levenberg-marquardt)}  \\ \cline{2-5}
  & \tabincell{c}{NJ (Iteration, time$/s$)}     & $||F(x^{it})||_\infty$  & \tabincell{c}{Iteration (time$/s$)} & $||F(x^{it})||_\infty$     \\\hline

  \tabincell{c}{1. Trid Function \cite{SB2013} \\ (m = 1999, n = 2000)}   & \tabincell{c}{2 (14, 1.2698$s$)} & 1.6764E-08 & \tabincell{c}{14 (87.1366$s$)} & 2.5611E-09   \\ \hline

  \tabincell{c}{2. Grewank Function \cite{SB2013} \\ (m = 1999, n=2000)}   & \tabincell{c}{12 (29, 9.7008$s$)} & 2.1524E-07 & \tabincell{c}{36 (308.7068$s$)} & 4.9506E-16 \\ \hline

  \tabincell{c}{3. Dixon Price Function \cite{SB2013} \\ (m = 1999, n = 2000)}   & \tabincell{c}{6 (21, 4.4584$s$)} & 8.6397E-09 & \tabincell{c}{7 (37.1861$s$)} & 1.5824E-08 \\ \hline

  \tabincell{c}{4. Rosenbrock Function \cite{MGH1981}  \\  (m = 1999, n = 2000)}   & \tabincell{c}{6 (22, 4.5985$s$)} & 4.8484E-09 & \tabincell{c}{7 (36.7351$s$)} & 1.2921E-09  \\ \hline

  \tabincell{c}{5. Trigonometric Function \cite{MGH1981}  \\  (m = 1999, n = 2000)}   & \tabincell{c}{13 (30, 3.2750$s$)} & 3.1687E-07 & \tabincell{c}{\red{70 (605.1742$s$)}} & \tabincell{c}{\red{2.2600e-05} \\ \red{(failed)}}  \\ \hline

  \tabincell{c}{6. Singular Broyden Function \cite{Luksan1994} \\ (m = 1999, n = 2000)}   & \tabincell{c}{19 (35, 13.0141$s$)} & 5.2082E-07  & \tabincell{c} {\red {12 (63.8643$s$)}}  & \tabincell{c} { \red{0.0213} \\ \red{(failed)}}\\ \hline

  \tabincell{c}{7. Extended Powell Singular Function \cite{MGH1981} \\  (m = 1999, n = 2000)}   & \tabincell{c}{10 (25, 7.3242$s$)} & 7.7570E-07 & \tabincell{c}{11 (60.1272$s$)} & 7.5422E-07  \\ \hline

  \tabincell{c}{8. Tridiagonal System Function \cite{Luksan1994} \\ (m = 1999, n = 2000)}   & \tabincell{c}{6 (20, 4.4418$s$)} & 5.3380E-07 & \tabincell{c}{8 (41.2931$s$)} & 1.6009E-14  \\ \hline

  \tabincell{c}{9. Discrete Boundary-Value Function \cite{Luksan1994} \\ (m = 1999, n = 2000)}   & \tabincell{c}{2 (13, 1.8657$s$)} & 1.8195E-07 & \tabincell{c}{6 (79.6152$s$)} & 1.5319E-12 \\ \hline

  \tabincell{c}{10. Broyden Tridiagonal Function \cite{Luksan1994} \\ (m = 1999, n = 2000)}   & \tabincell{c}{2 (15, 1.8283$s$)} & 3.0495E-07 & \tabincell{c}{5 (25.6165$s$)} & 1.4704E-10 \\ \hline

  \tabincell{c}{11. Extended Wood Function \cite{Andrei2008} \\ (m = 1999, n = 2000)}   & \tabincell{c}{6 (22, 4.5732$s$)} & 7.0474E-08 & \tabincell{c}{8 (41.4228$s$)} & 5.8341E-09  \\ \hline

  \tabincell{c}{12. Extended Cliff Function \cite{Andrei2008} \\ (m = 1999, n = 2000)}   & \tabincell{c}{6 (21, 4.6795$s$)} & 9.6747E-07 & \tabincell{c}{10 (98.8046$s$)} & 1.5507E-08 \\ \hline

  \tabincell{c}{13. Extended Hiebert Function \cite{Andrei2008} \\ (m = 1999, n = 2000)}   & \tabincell{c}{2 (16, 1.2882$s$)} & 1.4830E-08 & \tabincell{c}{2 (10.4330$s$)} & 8.0231E-10  \\ \hline

  \tabincell{c}{14. Extended Maratos Function \cite{Andrei2008} \\ (m = 1999, n = 2000)}   & \tabincell{c}{16 (37, 10.5885$s$)} & 7.4568E-07 & \tabincell{c} {\red {199 (1612.5726$s$)} } & \tabincell{c} {\red{0.1757}\\ \red{(failed)}} \\ \hline

  \tabincell{c}{15. Extended Psc1 Function \cite{Andrei2008} \\ (m = 1999, n = 2000)}   & \tabincell{c}{5 (22, 4.3262$s$)} & 2.0538E-07 & \tabincell{c}{7 (37.6980$s$)} & 8.4970E-13  \\ \hline

  \tabincell{c}{16. Extended Quadratic Penalty QP 1 \\ Function \cite{Andrei2008}  (m = 1999, n = 2000)}   & \tabincell{c}{10 (26, 7.4105$s$)} & 2.2681E-07 & \tabincell{c}{\red{11 (60.6542$s$)}} & \tabincell{c}{\red{1.5642E-04} \\ \red{(failed)}} \\ \hline

  \tabincell{c}{17. Extended Quadratic Penalty QP 2  \\ Function \cite{Andrei2008} (m = 1999, n = 2000)}   & \tabincell{c}{7 (32, 6.6672$s$)} & 4.6757E-08 & \tabincell{c}{\red{8 (42.7976$s$)}} & \tabincell{c}{\red{1.2390E-05} \\ \red{(failed)}} \\ \hline

  \tabincell{c}{18. Extended TET Function \cite{Andrei2008} \\ (m = 1999, n = 2000)}   & \tabincell{c}{7 (24, 5.9189$s$)} & 1.4731E-07 & \tabincell{c}{8 (41.0473$s$)} & 7.8525E-11 \\ \hline

  \tabincell{c}{19. EG2 Function \cite{Andrei2008} \\ (m = 1999, n = 2000)}   & \tabincell{c}{5 (22, 4.5872$s$)} & 3.1190E-08 & \tabincell{c}{153 (1225.0636$s$)} & 1.3816E-10 \\ \hline

  \tabincell{c}{20. Extended BD1 Function \cite{Andrei2008} \\ (m = 1999, n = 2000)}   & \tabincell{c}{5 (19, 3.9460$s$)} & 7.2238E-07 & \tabincell{c}{7 (38.3996$s$)} & 1.0549E-14 \\ \hline

\end{tabular}}
\end{table}

\vskip 2mm

\begin{table}[!http]
  \newcommand{\tabincell}[2]{\begin{tabular}{@{}#1@{}}#2\end{tabular}}
  \scriptsize
  \centering
  \caption{Numerical results of GCNMTr and fsolve for the determined problems with $m = 2000, \; n = 2000$.}
  \label{TABCOMT3}
  \resizebox{\textwidth}{!}{
  \begin{tabular}{|c|c|c|c|c|c|c|c|c|c|}
  \hline
  \multirow{2}{*}{Problems } & \multicolumn{2}{c|}{GCNMTr} & \multicolumn{2}{c|}{fsolve.m (levenberg-marquardt)}  \\ \cline{2-5}
  & \tabincell{c}{NJ (Iteration, time$/s$)}     & $||F(x^{it})||_\infty$  & \tabincell{c}{Iteration (time$/s$)} & $||F(x^{it})||_\infty$     \\\hline

  \tabincell{c}{1. Trid Function \cite{SB2013} \\ (m = 2000, n = 2000)}   & \tabincell{c}{2 (14, 1.1696$s$)} & 1.6764E-08 & \tabincell{c}{15 (84.5851$s$)} & 2.3283E-10   \\ \hline

  \tabincell{c}{2. Grewank Function \cite{SB2013} \\ (m = 2000, n=2000)}   & \tabincell{c}{37 (64, 26.9887$s$)} & 3.5969E-07 & \tabincell{c}{36 (299.4471$s$)} & 4.9853E-16 \\ \hline

  \tabincell{c}{3. Dixon Price Function \cite{SB2013} \\ (m = 2000, n = 2000)}   & \tabincell{c}{6 (21, 4.0241$s$)} & 8.6667E-09 & \tabincell{c}{7 (41.3923$s$)} & 1.5895E-08 \\ \hline

  \tabincell{c}{4. Rosenbrock Function \cite{MGH1981}  \\  (m = 2000, n = 2000)}   & \tabincell{c}{6 (22, 4.9494$s$)} & 6.5022E-08 & \tabincell{c}{ \red {199 (1671.5553$s$)}} & \tabincell{c}{\red{0.0043} \\ \red{(failed)}}  \\ \hline

  \tabincell{c}{5. Trigonometric Function \cite{MGH1981}  \\  (m = 2000, n = 2000)}   & \tabincell{c}{4 (13, 3.0073$s$)} & 3.4463E-07 & \tabincell{c}{\red{44 (420.8693$s$)}} & \tabincell{c}{\red{2.2750E-05} \\ \red{(failed)}}  \\ \hline

  \tabincell{c}{6. Singular Broyden Function \cite{Luksan1994} \\ (m = 2000, n = 2000)}   & \tabincell{c}{19 (35, 11.6583$s$)} & 5.5066E-07  & \tabincell{c} {\red {12 (75.0792$s$)}}  & \tabincell{c} { \red{0.0226} \\ \red{(failed)}}\\ \hline

  \tabincell{c}{7. Extended Powell Singular Function \cite{MGH1981}  \\  (m = 2000, n = 2000)}   & \tabincell{c}{10 (25, 6.9210$s$)} & 7.7570E-07 & \tabincell{c}{12 (73.0282$s$)} & 7.5422E-07  \\ \hline

  \tabincell{c}{8. Tridiagonal System Function \cite{Luksan1994} \\ (m = 2000, n = 2000)}   & \tabincell{c}{6 (20, 4.2114$s$)} & 1.7450E-07 & \tabincell{c}{7 (40.8813$s$)} & 1.8865E-11  \\ \hline

  \tabincell{c}{9. Discrete Boundary-Value Function \cite{Luksan1994} \\ (m = 2000, n = 2000)}   & \tabincell{c}{2 (13, 1.8843$s$)} & 6.8400E-07 & \tabincell{c}{13 (88.6350$s$)} & 1.0325E-11 \\ \hline

  \tabincell{c}{10. Broyden Tridiagonal Function \cite{Luksan1994} \\ (m = 2000, n = 2000)}   & \tabincell{c}{3 (15, 1.7556$s$)} & 8.1420E-07 & \tabincell{c}{5 (31.1220$s$)} & 1.3323E-15 \\ \hline

  \tabincell{c}{11. Extended Wood Function \cite{Andrei2008} \\ (m = 2000, n = 2000)}   & \tabincell{c}{6 (22, 4.2280$s$)} & 7.0474E-08 & \tabincell{c}{8 (49.8806$s$)} & 5.8341E-09  \\ \hline

  \tabincell{c}{12. Extended Cliff Function \cite{Andrei2008} \\ (m = 2000, n = 2000)}   & \tabincell{c}{6 (21, 4.4264$s$)} & 9.6747E-07 & \tabincell{c}{10 (102.9089$s$)} & 1.5507E-08 \\ \hline

  \tabincell{c}{13. Extended Hiebert Function \cite{Andrei2008} \\ (m = 2000, n = 2000)}   & \tabincell{c}{2 (16, 1.2165$s$)} & 1.4634E-08 & \tabincell{c}{2 (11.2343$s$)} & 5.0164E-12  \\ \hline

  \tabincell{c}{14. Extended Maratos Function \cite{Andrei2008} \\ (m = 2000, n = 2000)}   & \tabincell{c}{16 (37, 10.0705$s$)} & 7.4567E-07 & \tabincell{c} {\red {199 (1692.7199$s$)} } & \tabincell{c} {\red{0.1757}\\ \red{(failed)}} \\ \hline

  \tabincell{c}{15. Extended Psc1 Function \cite{Andrei2008} \\ (m = 2000, n = 2000)}   & \tabincell{c}{5 (20, 3.8124$s$)} & 6.1415E-07 & \tabincell{c}{7 (39.2547$s$)} & 2.4758E-14  \\ \hline

  \tabincell{c}{16. Extended Quadratic Penalty QP 1 \\ Function \cite{Andrei2008}  (m = 2000, n = 2000)}   & \tabincell{c}{11 (27, 6.6775$s$)} & 7.5538E-08 & \tabincell{c}{\red{11 (60.9623$s$)}} & \tabincell{c}{\red{5.2307E-04} \\ \red{(failed)}} \\ \hline

  \tabincell{c}{17. Extended Quadratic Penalty QP 2  \\ Function \cite{Andrei2008} (m = 2000, n = 2000)}   & \tabincell{c}{10 (32, 7.3550$s$)} & 3.4402E-09 & \tabincell{c}{\red{10 (115.7199$s$)}} & \tabincell{c}{\red{0.2043} \\ \red{(failed)}} \\ \hline

  \tabincell{c}{18. Extended TET Function \cite{Andrei2008} \\ (m = 2000, n = 2000)}   & \tabincell{c}{7 (22, 5.2577$s$)} & 2.1258E-07 & \tabincell{c}{8 (45.2060$s$)} & 7.8525E-11 \\ \hline

  \tabincell{c}{19. EG2 Function \cite{Andrei2008} \\ (m = 2000, n = 2000)}   & \tabincell{c}{5 (22, 4.2430$s$)} & 3.1007E-08 & \tabincell{c}{153 (1264.3533$s$)} & 1.2457E-10 \\ \hline

  \tabincell{c}{20. Extended BD1 Function \cite{Andrei2008} \\ (m = 2000, n = 2000)}   & \tabincell{c}{5 (19, 3.6184$s$)} & 1.1822E-07 & \tabincell{c}{6 (33.3613$s$)} & 1.4783E-08 \\ \hline

\end{tabular}}
\end{table}

\vskip 2mm

\vskip 2mm

\vskip 2mm

\vskip 2mm

\section{Conclusions}

In this article, we consider the generalized continuation Newton method with
the trust-region updating strategy and the new updating technique of the Jacobian
matrix for the underdetermined system (GCNMTr). For some large-scale
underdetermined and determined problems, numerical results show that GCNMTr is more robust
and faster than the traditional optimization method such as the
Levenberg-Marquardt method (the subroutine fsolve.m of the MATLAB R2020a
environment). The computational time of GCNMTr is about $1/8$ to $1/50$ of 
that of fsolve. We also analyze the global convergence and the local
superlinear convergence of the new method under the standard assumptions. From our
point of view, the generalized continuation Newton method (Algorithm \ref{ALGGCNMTR})
can be regarded as an alternative workhorse for the nonlinear equations and
we will extend it to the constrained nonlinear programming problems.

\vskip 2mm

\section*{Acknowledgments} The authors are grateful to two anonymous referees for 
their comments and suggestions which greatly improve presentation of this paper.  

\section * {Declarations}

\vskip 2mm

\noindent \textbf{Funding:} This work was supported in part by Grant
61876199 from National Natural Science Foundation of China, and Grant YJCB2011003HI
from the Innovation Research Program of Huawei Technologies Co., Ltd..

\vskip 2mm

\noindent \textbf{Conflicts of interest/Competing interests:} Not applicable.

\vskip 2mm

\noindent \textbf{Availability of data and material (data transparency):} If it is requested, we will
provide the test data.

\vskip 2mm

\noindent \textbf{Code availability (software application or custom code):} If it is requested, we will
provide the code.

\end{document}